\documentclass{article}

\newcommand{\bA}{{\mathbf{A}}}

\newcommand{\bB}{{\mathbf{B}}}

\newcommand{\bC}{{\mathbf{C}}}
\newcommand{\cC}{{\mathcal{C}}}

\newcommand{\cD}{{\mathcal{D}}}

\newcommand{\bfe}{{\mathbf{e}}}

\newcommand{\cE}{{\mathcal{E}}}

\newcommand{\bF}{{\mathbf{F}}}

\newcommand{\bG}{{\mathbf{G}}}
\newcommand{\cG}{{\mathcal{G}}}

\newcommand{\bH}{{\mathbf{H}}}

\newcommand{\bI}{{\mathbf{I}}}
\newcommand{\cI}{{\mathcal{I}}}

\newcommand{\bM}{{\mathbf{M}}}

\newcommand{\cN}{{\mathcal{N}}}

\newcommand{\cO}{{\mathcal{O}}}

\newcommand{\bP}{{\mathbf{P}}}

\newcommand{\bQ}{{\mathbf{Q}}}

\newcommand{\bR}{{\mathbf{R}}}
\newcommand{\cR}{{\mathcal{R}}}

\newcommand{\cS}{{\mathcal{S}}}

\newcommand{\bU}{{\mathbf{U}}}

\newcommand{\cU}{{\mathcal{U}}}

\newcommand{\bW}{{\mathbf{W}}}

\newcommand{\bX}{{\mathbf{X}}}

\newcommand{\bY}{{\mathbf{Y}}}

\newcommand{\bZ}{{\mathbf{Z}}}

\newcommand{\bpi}{{\boldsymbol{\pi}}}

\newcommand{\bTh}{{\boldsymbol{\Theta}}}

\newcommand{\bxi}{{\boldsymbol{\xi}}}

\newcommand{\bPi}{{\boldsymbol{\Pi}}}

\newcommand{\expect}{\mathbb{E}}
\newcommand{\reals}{\mathbb{R}}

\newcommand{\ip}[2]{\left\langle{#1},{#2}\right\rangle}

\newcommand{\avg}{\text{avg}}
\usepackage{arxiv}
\usepackage[applemac]{inputenc} 		
\usepackage[T1]{fontenc}    		
\usepackage[colorlinks]{hyperref}      
\usepackage{url}            			
\usepackage{amsthm}
\usepackage{booktabs}       		
\usepackage{amsfonts}       		
\usepackage{amsmath}
\usepackage{nicefrac}       		
\usepackage{microtype}      		
\usepackage{lipsum}				
\usepackage[square,numbers]{natbib}
\usepackage{mathtools}
\usepackage{algpseudocode}
\usepackage{graphicx}
\usepackage{indentfirst,latexsym,bm}
\usepackage{amsmath}
\usepackage{amssymb}
\usepackage{xcolor}
\usepackage{comment}
\usepackage{enumitem}
\usepackage{dsfont}
\usepackage{amssymb}

\usepackage{algorithm}
\usepackage{algorithmicx}
\usepackage{float}
\usepackage{subfig}
\usepackage{caption}
\usepackage{enumitem}
\usepackage{todonotes}
\usepackage{makecell}
\usepackage[flushleft]{threeparttable} 

\usepackage{tikz}
\usetikzlibrary{arrows.meta,positioning}
\usepackage{pgfplots}
\pgfplotsset{compat=newest}
\usepgfplotslibrary{fillbetween}
\usetikzlibrary{shapes,decorations}
\usetikzlibrary{fit}

\colorlet{color1}{blue}
\colorlet{color2}{red!50!black}

\definecolor{ivory}{RGB}{218,215,203}

\definecolor{cuhkp}{RGB}{98,56,105} 	
\definecolor{cuhkpl}{RGB}{152,24,147} 	
\definecolor{cuhkb}{RGB}{219,160,1} 	
\definecolor{cuhkbd}{RGB}{178,129,0} 	
\definecolor{cuhkr}{RGB}{88,35,155}  	

\usepackage{hyperref}[6.83]
\hypersetup{
  colorlinks=true,
  frenchlinks=false,
  pdfborder={0 0 0},
  naturalnames=false,
  hypertexnames=false,
  breaklinks,
  allcolors = cuhkbd,
  urlcolor = color2,
}

\RequirePackage[capitalize,nameinlink]{cleveref}[0.19]

\crefname{section}{section}{sections}
\crefname{subsection}{subsection}{subsections}

\Crefname{figure}{Figure}{Figures}

\crefformat{equation}{\textup{#2(#1)#3}}
\crefrangeformat{equation}{\textup{#3(#1)#4--#5(#2)#6}}
\crefmultiformat{equation}{\textup{#2(#1)#3}}{ and \textup{#2(#1)#3}}
{, \textup{#2(#1)#3}}{, and \textup{#2(#1)#3}}
\crefrangemultiformat{equation}{\textup{#3(#1)#4--#5(#2)#6}}%
{ and \textup{#3(#1)#4--#5(#2)#6}}{, \textup{#3(#1)#4--#5(#2)#6}}{, and \textup{#3(#1)#4--#5(#2)#6}}

\Crefformat{equation}{#2Equation~\textup{(#1)}#3}
\Crefrangeformat{equation}{Equations~\textup{#3(#1)#4--#5(#2)#6}}
\Crefmultiformat{equation}{Equations~\textup{#2(#1)#3}}{ and \textup{#2(#1)#3}}
{, \textup{#2(#1)#3}}{, and \textup{#2(#1)#3}}
\Crefrangemultiformat{equation}{Equations~\textup{#3(#1)#4--#5(#2)#6}}%
{ and \textup{#3(#1)#4--#5(#2)#6}}{, \textup{#3(#1)#4--#5(#2)#6}}{, and \textup{#3(#1)#4--#5(#2)#6}}

\crefdefaultlabelformat{#2\textup{#1}#3}

\theoremstyle{plain}
\newtheorem{theorem}{Theorem}[section]
\newtheorem{thm}{Theorem}[section]
\newtheorem{lemma}[thm]{Lemma}
\newtheorem{corollary}[thm]{Corollary}

\newtheorem{remark}{Remark}[section]
\newtheorem{assumption}[thm]{Assumption}
\theoremstyle{plain}

\newcommand{\prt}[1]{\left(#1\right)}
\newcommand{\brk}[1]{\left[#1\right]}
\newcommand{\crk}[1]{\left\{#1\right\}}

\newcommand{\norm}[1]{\left\|#1\right\|}

\usepackage{bbm}
\newcommand{\mone}{\mathbbm{1}}

\newcommand{\T}{\top}

\usepackage{etoolbox}
\AtBeginEnvironment{equation}{\small}
\AtBeginEnvironment{equation*}{\small}
\AtBeginEnvironment{displaymath}{\small}
\newcommand{\smallmath}[1]{\text{\small $#1$}}

\title{Distributed Learning over Arbitrary Topology: Linear Speed-Up with Polynomial Transient Time}

\author{
Runze You \\
School of Data Science (SDS) \\
The Chinese University of Hong Kong, Shenzhen \\
\texttt{runzeyou@link.cuhk.edu.cn}\\
\And
Shi Pu \\
School of Data Science (SDS) \\
The Chinese University of Hong Kong, Shenzhen \\
\texttt{pushi@cuhk.edu.cn}
}

\date{}

\begin{document}
\maketitle

\begin{abstract}
    We study a distributed learning problem in which $n$ agents, each with potentially heterogeneous local data, collaboratively minimize the sum of their local cost functions via peer-to-peer communication. We propose a novel algorithm, \emph{Spanning Tree Push-Pull} (STPP), which employs two spanning trees extracted from a general communication graph to distribute both model parameters and stochastic gradients. Unlike prior approaches that rely heavily on spectral gap properties, STPP leverages a more flexible topological characterization, enabling robust information flow and efficient updates. Theoretically, we prove that STPP achieves linear speedup and polynomial transient iteration complexity---up to $\cO(n^7)$ for smooth nonconvex objectives and $\tilde{\cO}(n^3)$ for smooth strongly convex objectives---under arbitrary network topologies. Moreover, compared with existing methods, STPP achieves faster convergence rates on sparse and non-regular topologies (e.g., directed rings) and reduces communication overhead on dense networks (e.g., static exponential graphs). 
    Numerical experiments further demonstrate the strong performance of STPP 
    across various graph architectures.
\end{abstract}

\section{Introduction}

Consider a set of networked agents, labeled as $\cN = \{1,2,\dots,n\}$, where each agent~$i$ has a local cost function $f_i:\reals^p \rightarrow \reals$ and can only communicate with its immediate neighbors. These agents jointly seek to minimize the average of their local cost functions
\begin{equation}
    \label{eq:obj}
    \min_{x \in \reals^p} \; f(x) 
    \;:=\;
    \frac{1}{n}\sum_{i=1}^n f_i(x),
\end{equation}
where $f_i(x) := \expect_{\xi_i\sim \cD_i} F_i(x;\xi_i)$. Here, $\xi_i$ denotes the local data following the distribution $\cD_i$. \emph{Data heterogeneity} arises if $\{\cD_i\}_{i=1}^n$ are not identical. The above formulation also covers \emph{empirical risk minimization} (ERM), with each $\cD_i$ representing a local training dataset.

To solve \eqref{eq:obj}, we assume each agent can query a stochastic oracle ($\cS\cO$) to obtain noisy gradient samples. Such \emph{stochastic gradients} appear widely in online distributed learning \cite{recht2011hogwild}, reinforcement learning \cite{mnih2013playing, lillicrap2015continuous}, generative modeling \cite{goodfellow2014generative, kingma2019introduction}, and parameter estimation \cite{bottou2012stochastic, srivastava2014dropout}. Stochastic Gradient Descent (SGD) is commonly used in modern machine learning, where problems often involve massive datasets and high-dimensional models. However, standard SGD operates sequentially, making it difficult to scale to very large datasets~\cite{dean2012large}. Consequently, a variety of \emph{distributed} algorithms for \eqref{eq:obj} have been developed to improve scalability; see, e.g., \cite{pu2021distributed,lian2017can,ding2023dsgd,ying2021exponential}.

A  classical approach for distributed learning is the \emph{centralized} master-worker paradigm, in which each worker node communicates with a (possibly virtual) central server \cite{zinkevich2010parallelized,li2014scaling}. However, this often leads to significant communication overheads and high latency due to limited bandwidth between the server and workers \cite{assran2019stochastic,nadiradze2021asynchronous}.

\emph{Decentralized learning} has emerged as an attractive alternative for reducing communication costs. In a decentralized setup, the nodes form a specific network topology (e.g., ring, grid, exponential graph) and exchange information \emph{locally} with only their direct neighbors \cite{assran2019stochastic}. Such communication patterns greatly reduce overhead, and they also accommodate scenarios in which certain nodes can only reach a subset of the network (due to varying power ranges or physical restrictions)~\cite{yang2019survey}.

The critical factor that distinguishes decentralized from centralized methods is the \emph{communication graph} (or network topology).  The graph structure strongly influences the convergence rate: dense topologies (e.g., exponential graphs) generally yield faster convergence compared with sparse topologies (e.g., rings) \cite{pu2021sharp,nedic2018network}. Importantly, for the stochastic  gradient setting, it has been shown that decentralized algorithms can eventually match the convergence rate of centralized methods, but only after incurring a period of so-called \emph{transient iterations}~\cite{lian2017can,ying2021exponential}. These transient iterations heavily depend on the underlying network topology, making them a crucial measure of an algorithm’s efficiency.

In this work, we explore an alternative to standard gossip-based methods by leveraging a \emph{spanning tree} construction, inspired by \cite{pu2020push}. Rather than using a single graph for message passing at every iteration, our proposed \emph{Spanning Tree Push-Pull} (STPP) algorithm employs two spanning trees: a \emph{Pull Tree} $\cG_{\bR}$ and a \emph{Push Tree} $\cG_{\bC}$. Model parameters are transmitted through $\cG_{\bR}$ from parent nodes to child nodes, while the associated stochastic gradients (computed using the current parameters) are \emph{accumulated} through $\cG_{\bC}$. Conceptually, each agent in the spanning trees acts like a worker on an assembly line, passing model updates and gradients along separate paths. 

STPP can be viewed as \emph{semi-(de)centralized} due to the pivotal role of node~1 in both spanning trees. Notably, the corresponding mixing matrices for $\cG_{\bR}$ and $\cG_{\bC}$ are composed solely of 0’s and 1’s, which simplifies updates and improves efficiency. We show that STPP attains \emph{linear speedup} and \emph{polynomial transient iteration complexity}—up to $\cO(n^7)$ for smooth nonconvex objectives and $\tilde{\cO}(n^3)$ for smooth strongly convex objectives—across arbitrary network topologies. Furthermore, STPP outperforms existing algorithms on several commonly used topologies, as evidenced by the results in Table \ref{tab:smooth} and Table \ref{tab:convex}.

\subsection{Related Works}

\textbf{Decentralized Learning.}
Decentralized Stochastic Gradient Descent (DSGD) and its variants \cite{lian2017can,ying2021exponential,koloskova2019decentralized} have gained significant attention for large-scale distributed training. While these methods are flexible and often effective, they can suffer under data heterogeneity \cite{koloskova2020unified}, necessitating more advanced techniques such as EXTRA \cite{shi2015extra}, Exact-Diffusion/$\text{D}^2$ \cite{li2019decentralized}, and gradient tracking \cite{pu2021distributed}. 
A broad family of Push-Pull methods \cite{pu2020push,xin2018linear,xin2020general} has been proposed to further relax topological requirements, albeit initially for deterministic settings and strongly convex objectives. 
In this work, we build upon the flexibility of Push-Pull by exploiting \emph{spanning trees} in the underlying network.

\textbf{Network Topology and Weight Matrices.}
Most decentralized learning algorithms rely on specific constraints for mixing matrices (e.g., doubly stochastic) that may be challenging or impossible to satisfy in certain directed graphs \cite{gharesifard2010does}.  
Push-Sum (SGP) \cite{assran2019stochastic} and Push-DIGing \cite{xi2017add,nedic2017achieving,liang2023towards} techniques adopt column-stochastic mixing, enabling consensus on directed strongly connected graphs, but can incur prohibitive transient iteration complexities (e.g., $\cO(2^n)$ for some dense networks \cite{liang2023towards}).  
Some alternative approaches, such as \cite{zhang2019fully,tian2020achieving}, relax the requirements for mixing matrices by employing two distinct weight matrices through gradient tracking, which still require the communication graph to be strongly connected. 
Recent Push-Pull methods \cite{pu2020push,xin2018linear,zhao2023asymptotic} operate on two separate subgraphs that each contain a spanning tree and share a common root.  BTPP \cite{you2024b} uses exactly two spanning trees from the B-ary Tree family. However, existing results either lack clear establishment of linear speedup or fail to provide tight complexity bounds applicable to arbitrary network topologies.   
In contrast, STPP leverages spanning trees to provide more favorable transient iteration complexity and broader applicability, as shown in Table \ref{tab:smooth} (smooth nonconvex objectives) and Table \ref{tab:convex} (smooth strongly convex objectives).

\begin{table*}[ht]
    \caption{Comparison of  distributed stochastic gradient methods by their transient iteration bounds under smooth nonconvex objectives. The notation $\tilde{\cO}(\cdot)$ hides polylogarithmic factors. ``N.P.'' indicates non-polynomial bounds. We list the best-known ``DSGT/Push-DIGing'' results from \cite{alghunaim2022unified,nedic2018network,koloskova2020unified,pu2021distributed,liang2023towards} and ``DSGD/SGP'' results from \cite{ying2021exponential,lian2017can,assran2019stochastic}. ``General Topo.'' refers to arbitrary network topologies.}
    \label{tab:smooth}
    \begin{center}
    \begin{small}
    \begin{sc}
    \begin{tabular}{cccccc}
    \toprule
    Algorithm & General Topo. & Di. Ring & Ring & Grid  & Static Exp.  \\
    \midrule
    DSGD/SGP & N.P. & $\cO(n^{11})$ & $\cO(n^{11})$ & $\tilde{\cO}(n^{7})$   & $\tilde{\cO}(n^3)$  \\
    DSGT/Push-DIGing & N.P. & $\cO(n^{11})$ & $\cO(n^{7})$ &  $\tilde{\cO}(n^{5})$  & $\tilde{\cO}(n^3)$  \\
    \textbf{STPP (Ours)} & $\boldsymbol{\cO(n^7)}$ & $\boldsymbol{\cO(n^7)}$ & $\boldsymbol{\cO(n^7)}$ & $\boldsymbol{\cO(n^4)}$  & $\boldsymbol{\tilde{\cO}(n)}$   \\
    \bottomrule
    \end{tabular}
    \end{sc}
    \end{small}
    \end{center}
\end{table*}
    
\begin{table*}[ht]
    \caption{Comparison of  distributed stochastic gradient methods under smooth strongly convex objectives. ``N.A.'' indicates cases not reported in the corresponding references.  We list the best-known results from \cite{pu2021sharp,song2022communication, koloskova2020unified,alghunaim2022unified,nedic2018network}.}
    \label{tab:convex}
    \begin{center}
    \begin{small}
    \begin{sc}
    \begin{tabular}{cccccc}
    \toprule
    Algorithm & General Topo. & Di. Ring & Ring & Grid  & Static Exp. \\
    \midrule
    DSGD/SGP & N.A. & $\tilde{\cO}(n^{5})$ & $\tilde{\cO}(n^{5})$ & $\tilde{\cO}(n^{3})$ &$\tilde{\cO}(n)$ \\
    DSGT/Push-DIGing & N.A. & $\tilde{\cO}(n^{5})$ & $\tilde{\cO}(n^{3})$ &  $\tilde{\cO}(n^{2})$   & $\tilde{\cO}(n)$ \\
    \textbf{STPP (Ours)} & $\boldsymbol{\tilde{\cO}(n^3)}$ & $\boldsymbol{\tilde{\cO}(n^3)}$ & $\boldsymbol{\tilde{\cO}(n^3)}$ & $\boldsymbol{\tilde{\cO}(n^{3/2})}$  & $\boldsymbol{\tilde{\cO}(1)}$   \\
    \bottomrule
    \end{tabular}
    \end{sc}
    \end{small}
    \end{center}
\end{table*}

\subsection{Main Contribution}
In this paper, we introduce a novel distributed learning algorithm, called \emph{spanning Tree Push-Pull} (\textbf{STPP}), which effectively leverages spanning trees extracted from the overall communication graph. Our approach remains provably efficient for solving Problem~\eqref{eq:obj} under general network topologies with minimal connectivity requirements. The main contribution of this work is summarized as follows:
\begin{itemize}
    \item \textbf{New Topology Characterization.}
    We show that STPP achieves linear speedup for arbitrary network topologies and outperforms prior methods under common topologies. Unlike existing analyses that often rely solely on spectral gaps \cite{alghunaim2022unified,koloskova2021improved,nedic2018network}, we demonstrate that two new graph characteristics---the average distance to the root node and the diameter of the spanning trees---play a critical role in algorithmic performance. This insight into the network structure provides a more comprehensive understanding of the convergence behavior.

    \item \textbf{Polynomial Transient Iteration Complexity.}
    We prove that STPP enjoys a polynomial transient iteration complexity under arbitrary network topologies. Specifically, for smooth nonconvex objectives, the worst-case complexity is upper bounded by $\mathcal{O}(n^7)$, while for strongly convex objectives it is $\tilde{\mathcal{O}}(n^3)$. These bounds significantly outperform those of competing algorithms (see  Tables \ref{tab:smooth} and \ref{tab:convex}). For large $n$, the transient iteration complexity has a notable impact on real-world performance, as confirmed by our numerical experiments.

    \item \textbf{Novel Analytical Techniques.}
    The convergence analysis of STPP poses unique challenges because the algorithm employs two distinct topologies for model parameters and gradient trackers. Rather than constructing induced matrix norms $\|\cdot\|_{\mathbf{R}}$ and $\|\cdot\|_{\mathbf{C}}$ as in \cite{pu2020push,zhu2024r,zhao2023asymptotic}, we carry out the analysis under standard Euclidean norms ($\|\cdot\|_{2}$ and $\|\cdot\|_{F}$). Through carefully managing matrix products and associated terms, we establish theoretical guarantees that highlight the robust and efficient nature of STPP.
\end{itemize}

\subsection{Notations}
Throughout the paper, vectors default to columns if not otherwise specified. Let each agent $i$ hold a local copy $x_i \in \reals^p$ of the decision variable and an auxiliary variable $y_i\in \reals^p$. Their values at iteration $k$ are denoted by $x_i^{(k)}$ and $y_i^{(k)}$, respectively. We let $\bX = \brk{x_1, x_2, \cdots,x_n}^\T \in \reals^{n\times p},\bY = \brk{y_1, y_2, \cdots,y_n}^\T \in \reals^{n\times p}$ and $\mone$ denotes th column vector with all entries equal to 1. In addition, denote $\bxi := \brk{\xi_1, \xi_2, \cdots, \xi_n}^\T$ and $\bG(\bX,\bxi):= \brk{g_1(x_1,\xi_1), g_2(x_2,\xi_2), \cdots, g_n(x_n,\xi_n)}^\T \in \reals^{n\times p}$ as the aggregations of random samples and stochastic  gradients respectively. To make the proof more concise, we use $\bG^{(t)}$ to represent $\bG(\bX^{(t)},\bxi^{(t)})$, and define an aggregate gradients of the local variables as 
$
\begin{aligned}
    \nabla F(\bX) :=  \brk{  \nabla f_1 (x_1),\nabla f_2 (x_2),  \cdots, \nabla f_n(x_n) }^\T \in \reals^{n\times p} ,
\end{aligned}
$
where $\bF(\bX):= \sum_{i=1}^{n} f_i(x_i)$. Besides, we use $\nabla\bF^{(t)}$ to represent $\nabla \bF(\bX^{(t)})$ and $\bTh^{(t)}:=  \bG^{(t)} - \nabla\bF^{(t)}$ to represent the gradient estimation error.  The term $\ip{a}{b}$ stands for the inner product of two vectors $a,b\in \reals^p$. For matrices, $\norm{\cdot}_2$ and $\norm{\cdot}_F$ represent the spectral norm and the Frobenius norm respectively, which degenerate to the Euclidean norm $\norm{\cdot}$ for vectors. For simplicity, any non-zero square matrix with power $0$ is the unit matrix $\bI$ with the same dimension if not otherwise specified.

A directed graph $\cG(\cN, \cE)$ consists of a set of $n$ nodes $\cN$ and a set of directed edges $\cE\subseteq \cN\times \cN$, where an edge $(j,i)\in \cE$ indicates that node $j$ can directly send information to node $i$. To facilitate the local averaging procedure, each graph can be associated with a non-negative weight matrix $\bW = \brk{w_{ij}} \in \reals^{n\times n}$, whose element $w_{ij}$ is non-zero only if $(j, i)\in \cE$. 
Similarly, a non-negative weight matrix $\bW$ corresponds to a directed graph denoted by $\cG_\bW$.
For a given graph $\cG_\bW$, the in-neighborhood and out-neighborhood of node $i\in \cN$ are given by $\cN_{\bW,i}^{\text{in}}:= \{j\in\cN: (j,i)\in\cE\}$ and $\cN_{\bW,i}^{\text{out}}:= \{j\in\cN: (i,j)\in\cE\}$, respectively. 

\subsection{Organization of the Paper}
In Section \ref{sec:algo}, we present the proposed Spanning Tree Push-Pull  (STPP) algorithm. In Section \ref{sec:convergence}, we provide the main convergence results of the algorithm. In Section \ref{sec:exp}, we demonstrate the effectiveness of the proposed STPP algorithm through numerical experiments. Finally, we conclude the paper in Section \ref{sec:conclusion}.

\section{Spanning Tree Push-Pull Method}
\label{sec:algo}

 In this section, we introduce the Spanning Tree Push-Pull  (STPP) algorithm with special focus on the communication graphs and the properties of the mixing matrices.
\subsection{Communication Graphs}

We consider a cooperative setting where all agents aim to solve Problem \eqref{eq:obj} over a network modeled as a directed graph $\cG(\cN,\cE)$, where $\cN = \{1, 2, \dots, n\}$ denotes the set of nodes and $\cE \subset \cN \times \cN$ denotes the set of edges (communication links).
\begin{assumption}
    \label{a.graph}
    The communication graph $\cG(\cN, \cE)$ is strongly connected, i.e., there exists a directed path between every pair of nodes in the graph.
\end{assumption}
 Assumption \ref{a.graph} has been employed in most existing works \cite{liang2023towards, pu2021sharp, alghunaim2022unified, pu2020push, zhu2024r, xin2018linear}. To improve the flexibility in the topology design, given a graph $\cG(\cN, \cE)$ satisfying Assumption \ref{a.graph}, we can construct two spanning trees: the \emph{Pull Tree} $\cG_{\bR}$ and the \emph{Push Tree} $\cG_{\bC}$, by removing certain edges from $\cG(\cN, \cE)$. Both $\cG_{\bR}$ and $\cG_{\bC}$ share the same designated root node, labeled as node~1 for convenience. The existence of these spanning trees is guaranteed by Assumption \ref{a.graph}, as there exists a path between any pair of nodes in both directions.

In the Pull Tree $\cG_{\bR}$, each node (except the root) has exactly one parent node and potentially several child nodes. Node~1, the root, has no parent. In the Push Tree $\cG_{\bC}$, each node (again, except the root) has exactly one child node and possibly multiple parent nodes. Note that $\cG_{\bC}$ need \emph{not} simply be the edge-reversed version of $\cG_{\bR}$.

\begin{figure}[htbp]
	\centering
	\subfloat{\includegraphics[width=\linewidth]{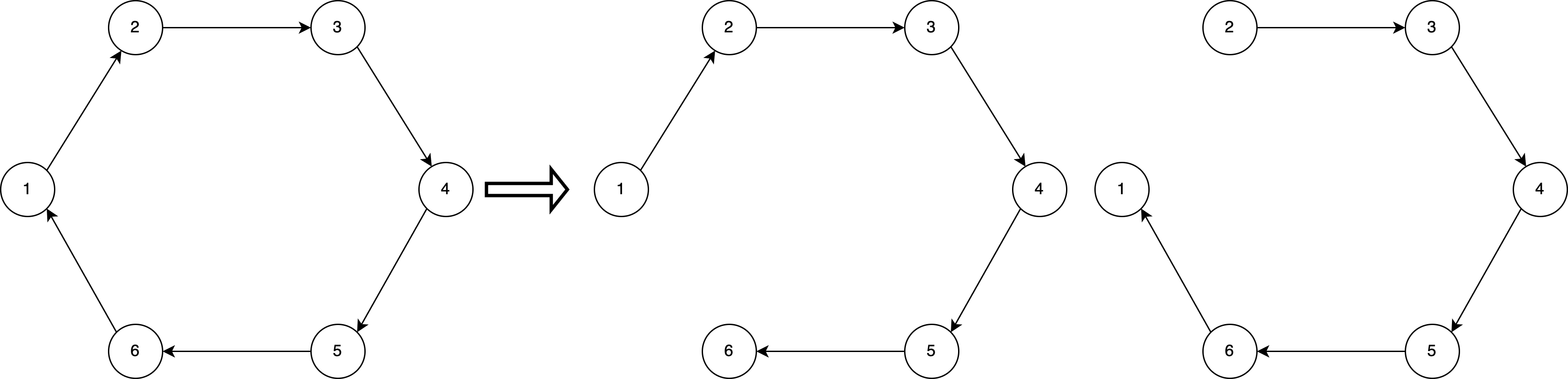}}
    \hfill
	\subfloat{\includegraphics[width=\linewidth]{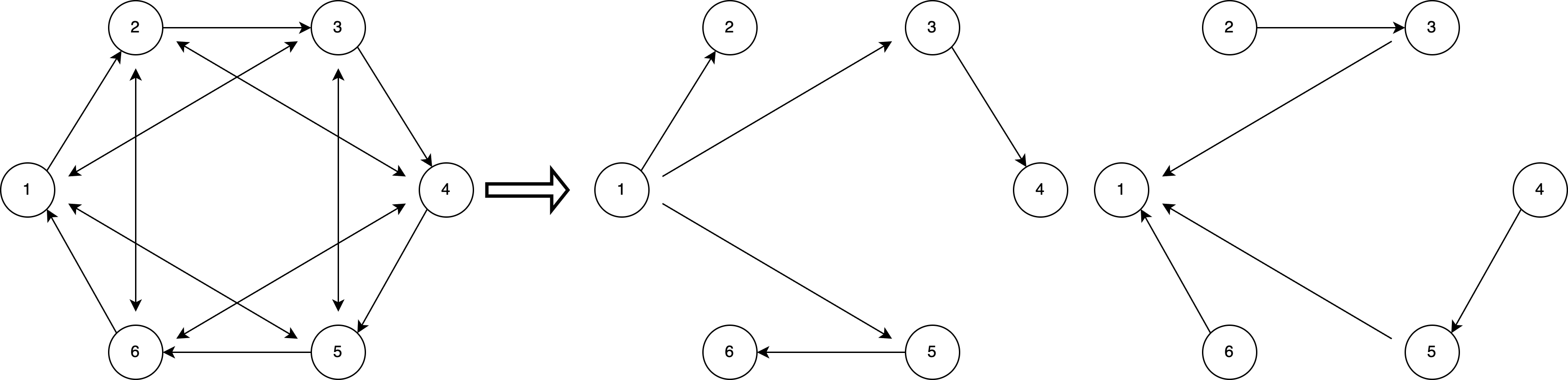}}
	\caption{\label{fig:ring3} 
    Illustration of network structures with 6 nodes. 
    \textbf{Top}: Directed ring graph (left) and its derived Pull Tree $\cG_{\bR}$ (center) and Push Tree $\cG_{\bC}$ (right), rooted at node 1. \textbf{Bottom}: Exponential graph (left) with corresponding Pull Tree $\cG_{\bR}$ (center) and Push Tree $\cG_{\bC}$ (right), also rooted at node 1.
    }
\end{figure}

The left panels of Figure~\ref{fig:ring3} depict two directed graphs (ring and exponential) with $n=6$ nodes. The two panels on the right illustrate how selectively removing edges from these graphs yields the corresponding Pull Tree $\cG_{\bR}$ and Push Tree $\cG_{\bC}$, both rooted at node~1.

To further characterize these spanning trees,  we define the \emph{distance} between two nodes $i$ and $j$ in the tree as the number of edges in the shortest path between them (length~1 if $i$ is the direct child of $j$). The \emph{diameter} of $\cG_{\bR}$ or $\cG_{\bC}$ is the maximum distance from any node to the root node~1. We denote these diameters by $d_\bR$ and $d_\bC$, respectively.

\subsection{Algorithm}

The proposed \emph{Spanning Tree Push-Pull} (STPP) algorithm for solving Problem~\eqref{eq:obj} is presented in Algorithm \ref{alg:stpp}, implemented on the communication graph $\cG(\cN, \cE)$. 

STPP initializes with two sub-graphs, the Pull Tree $\cG_{\bR}$ and the Push Tree $\cG_{\bC}$, from the communication graph $\cG(\cN, \cE)$, which are used to facilitate the information exchange among agents.  Let the common root node of $\cG_\bR$ and $\cG_{\bC}$ be node $1$.

At each iteration $t$, each agent $i$ pulls the state information from its in-neighborhood $\cN_{\bR,i}^{\text{in}}$, pushes its gradient tracker to the out-neighborhood $\cN_{\bC,i}^{\text{out}}$, and updates its local variables $x_i$ and $y_i$ based on the received information. The agents aim to find the $\epsilon$-stationary point  or optimal solution jointly by performing  local computation and information exchange through $\cG_{\bR}$ and $\cG_{\bC}$.

\begin{algorithm}[htbp]
   \caption{\textbf{S}panning \textbf{T}ree \textbf{P}ush-\textbf{P}ull Algorithm (\textbf{STPP})}
   \label{alg:stpp}
\begin{algorithmic}
    \State {\bfseries Input:}  Communication graphs $\cG_\bR$ (Pull Tree) and $\cG_\bC$ (Push Tree). Each agent $i$ initializes with arbitrary but identical $x_{i}^{(0)} = x^{(0)} \in \reals^p$, stepsize $\gamma > 0$, initial gradient tracker $y_i^{(0)} = g_i(x_i^{(0)}, \xi_i^{(0)})$ by drawing a random sample $\xi_i^{(0)}$.
    \For{Iteration $t = 0$ {\bfseries to} $T-1$}
    \For{Agent $i$ {\bfseries in parallel}}
    \State  Pull $x_j^{(t)} - \gamma y_{j}^{(t)}$ from each $j \in \cN_{\bR,i}^{\text{in}}$.
    \State Push $y_{i}^{(t)}$ to each $j\in \cN_{\bC,i}^{\text{out}}$.
    \State Independently draw a random sample $\xi_i^{(t+1)}$.
    \State \textbf{Update} parameters through
        \[
        \begin{aligned}
            x_i^{(t+1)} & = \sum_{j \in \cN_{\bR,i}^{\text{in}} } \prt{x_j^{(t)} - \gamma y_{j}^{(t)}}, \\
            y_i^{(t+1)} & = \sum_{j\in\cN_{\bC,i}^{\text{in}} } y_j^{(t)} +  g_i(x^{(t+1)}_i;\xi_i^{(t+1)})  -  g_i(x^{(t)}_i;\xi_i^{(t)}). 
        \end{aligned}
        \]
    \EndFor
    \EndFor
\State {\bfseries Output:} $x_{1}^{(T)}$.
\end{algorithmic}
\end{algorithm}

\subsection{Weight Matrices}

We now construct the weight matrices $\bR$ and $\bC$ that encode the  topologies of $\cG_{\bR}$ and $\cG_{\bC}$, respectively. First, let $\cI_{\bR, i}^{k}$ be the set of nodes at distance~$k$ from node~$i$ in $\cG_{\bR}$ (i.e., $k$ edges away). By convention, set $\cI_{\bR, i}^{0} = \{i\}$. A similar definition applies for $\bC$ in the Push Tree.
Using these sets, define the matrices $\bR$ and $\bC$ as follows:
\begin{equation}
    \label{eq:rc}
    \begin{aligned}
        \bR_{i,j} & = 
        \begin{cases}
        1 & \text{if } i \in \cI_{\bR,j}^{1} \text{ or } (i=j=1), \\
        0 & \text{otherwise},
        \end{cases}
        \\
        \bC_{i,j} & =
        \begin{cases}
        1 & \text{if } i \in \cI_{\bC,j}^{1} \text{ or } (i=j=1), \\
        0 & \text{otherwise}.
        \end{cases}
    \end{aligned}
\end{equation}

Hence, $\bR$ is row-stochastic while $\bC$ is column-stochastic, each comprised solely of 0's and 1's. For instance, the mixing matrices corresponding to the trees derived from the directed ring graph in Figure \ref{fig:ring3} have the forms
\[
\bR = \prt{\begin{matrix}
    1 &   &   &   &   &   \\
    1 &   &   &   &   &   \\
      & 1 &   &   &   &   \\
      &   & 1 &   &   &   \\
      &   &   & 1 &   &   \\
      &   &   &   & 1 &   \\
    \end{matrix} },  \bC = \prt{\begin{matrix}
        1 &   &   &   &   & 1 \\
          &   &   &   &   &   \\
          & 1 &   &   &   &   \\
          &   & 1 &   &   &   \\
          &   &   & 1 &   &   \\
          &   &   &   & 1 &   \\
    \end{matrix} }, 
\]
where unspecified entries are zero.

\textbf{Closed-form analysis of powers of \(\bR\).} Let $\bfe_{\cI_{\bR,i}^{k}} \in \reals^n$ be a \emph{column indicator vector} taking the value 1 in the rows indexed by $\cI_{\bR,i}^{k}$ and 0 elsewhere. It follows that the $j$-th column of $\bR$ is precisely $\bfe_{\cI_{\bR,j}^{1}}$. Define
\[
\tilde{\cI}_{\bR,i}^{k} := \cup_{m\le k} \cI_{\bR,i}^{m},
\]
i.e., the set of all nodes within distance at most~$k$ from node~$i$. Construct a matrix $\bZ_k \in \reals^{n\times n}$ by arranging these indicator vectors as columns:
\[
\bZ_k := \brk{\bfe_{\tilde{\cI}_{\bR,1}^{k}}, \bfe_{\cI_{\bR,2}^{k}}, \bfe_{\cI_{\bR,3}^{k}}, \cdots, \bfe_{\cI_{\bR,n}^{k}}}.
\]
Lemma \ref{lem:Rk} shows that $\bR^k$ admits a simple expression in terms of $\bZ_k$.
\begin{lemma}
    \label{lem:Rk}
    For the pull matrix $\bR$ associated with $\cG_{\bR}$, and any positive integer $k$, we have
    \[
    \bR^k = \bZ_k.
    \]
\end{lemma}
\begin{proof}
   See Appendix \ref{pf:lem:Rk}.
\end{proof}

\textbf{Eigenvectors and diameters.}
Let $\bpi_{\bR}$ be the (normalized) left eigenvector of $\bR$ (associated with eigenvalue 1), and let $\bpi_{\bC}$ be the (normalized) right eigenvector of $\bC$ (also associated with eigenvalue 1). By Lemma \ref{lem:Rk}, we have
\[
\bpi_{\bR} \;=\;\bpi_{\bC} \;=\; \bigl(1,0,\dots,0\bigr)^\T \in \reals^n.
\]
Lemma \ref{lem:Rk_d} shows that once $k$ exceeds the diameter $d_{\bR}$ (or $d_{\bC}$ for the push tree), $\bR^k$ (resp., $\bC^k$) collapses into a matrix broadcasting (resp., converging) only from (resp., to) node~1.
\begin{corollary}
\label{lem:Rk_d}
Let $d_\bR$ and $d_\bC$ be the diameters of $\cG_{\bR}$ and $\cG_{\bC}$, respectively. Then
\[
\bR^{d_\bR} =\mone\bpi_{\bR}^\T,
\text{\, and\, }
\bC^{d_\bC} = \bpi_{\bC}\mone^\T.
\]
\end{corollary}
\begin{proof}
    By the definition of the diameter, we have $\tilde{\cI}_{\bR,1}^{d_\bR} = \crk{1,2,\cdots,n}$, which implies that $\bfe_{\tilde{\cI}_{\bR,1}^{d_{\bR}}} = \mone$. Considering that $\bR$ is row stochastic and nonnegative, i.e., $\bR\mone = 1$, we have $\bfe_{\cI_{\bR,j}^{d_{\bR}}} = \mathbf{0}$ holds for $j > 1$. Thus, $\bR^{d_\bR} = \mone \bpi_{\bR}^\T$. We also have $\bC^{d_\bC} = \bpi_{\bC} \mone^\T$ as $\bC^\T$ shares similar property as $\bR$.
\end{proof}

\begin{figure}[ht]
    \begin{center}
    \centerline{\includegraphics[width=0.9\columnwidth]{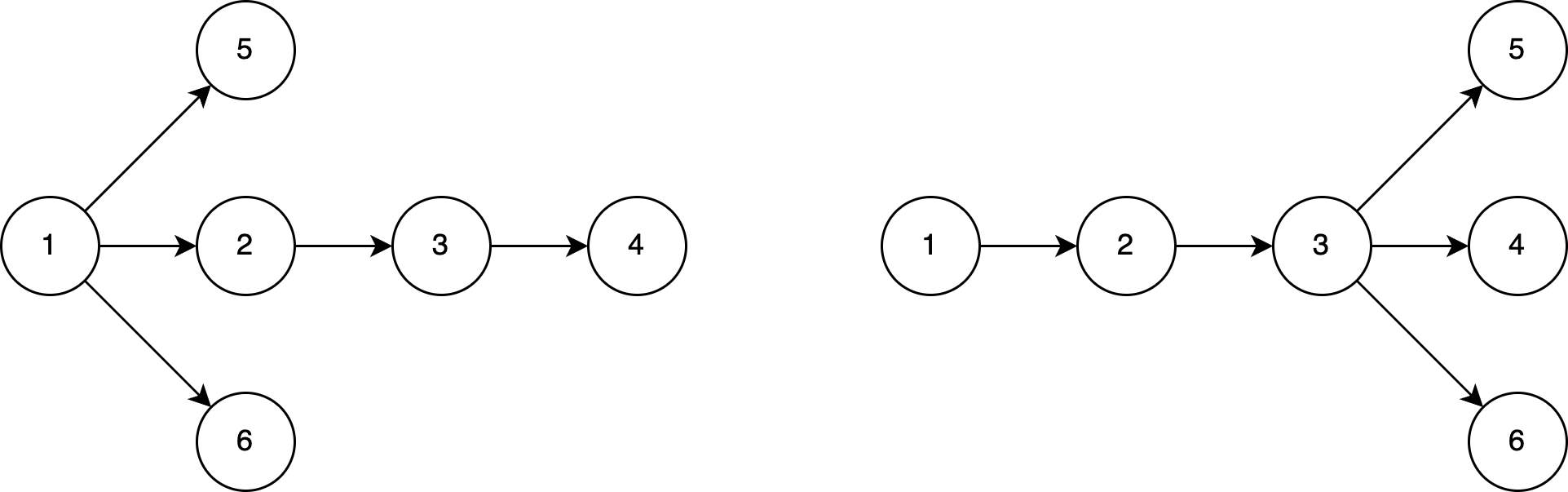}}
    \caption{Two Pull Trees ($n=6$) with the same diameter but different average distances from node~1.}
    \label{fig:r_tree}
    \end{center}
\end{figure}

\textbf{Average distances and bounded norms.}
Although the graph diameter captures the smallest distance to the root, it does not fully capture the performance difference when two trees share the same diameter but have different ``shapes''. Figure \ref{fig:r_tree} provides an example with $n=6$ where both spanning trees have the same diameter but differ in average distances from node~1.
We thus define the average distances $r_{\avg}$ and $c_{\avg}$ for $\cG_{\bR}$ and $\cG_{\bC}$, respectively, via
\[
r_k =\mone^\T\bfe_{\tilde{\cI}_{\bR,1}^k},
\quad
c_k =\mone^\T\bfe_{\tilde{\cI}_{\bC,1}^k},
\]
where $r_k$ (resp.\ $c_k$) is the number of nodes within distance $<k$ of node~1 in $\cG_{\bR}$ (resp.\ $\cG_{\bC}$). We set $r_0=c_0=1$. Then,
\[
r_{\avg}
=
\frac{1}{n}\sum_{k=0}^{d_\bR-1}\!\bigl(n-r_k\bigr),
\quad
c_{\avg}
=
\frac{1}{n}\sum_{k=0}^{d_\bC-1}\!\bigl(n-c_k\bigr).
\]
Lemma \ref{lem:Rk_norm} shows how these quantities bound the $2$-norm of $\bR^k$ and $\bC^k$ before $k$ reaches the diameters.
\begin{lemma}
    \label{lem:Rk_norm}
    Given the power $k$, we have
    \[
    \norm{\bR^k - \mone\bpi_{\bR}^\T}_2^2 \left\{ \begin{aligned}
        \le &\quad 2\prt{n-r_{k}} & \text{ if } k \le d_\bR -1, \\
        = &\quad 0 & \text{ if } k \ge  d_\bR,
    \end{aligned}   \right. \quad 
    \]
    and 
    \[
    \norm{\bC^k - \bpi_{\bC}\mone^\T}_2^2 \left\{ \begin{aligned}
        \le & \quad 2\prt{n-c_{k}} & \text{ if } k \le d_\bC -1, \\
        = &\quad  0 & \text{ if } k \ge  d_\bC .
    \end{aligned}   \right. \quad 
    \]
\end{lemma}
\begin{proof}
    By using the inequality $\norm{\cdot}_2\le \norm{\cdot}_F$ and applying Lemma \ref{lem:Rk}, we have, for any $k\le d_\bR - 1$, 
    \[
    \begin{aligned}
         & \norm{\bR^k - \mone\bpi_{\bR}^\T}_2^2 \le \norm{\bR^k - \mone\bpi_{\bR}^\T}_F^2 =  2 \norm{\bfe_{\tilde{\cI}_{\bR,1}^{k}} - \mone}^2_2 = 2\prt{n - r_k}.
    \end{aligned}
    \]
    Similar results hold for the term $\bC^k - \bpi_{\bC}\mone^\T$.
\end{proof}

\textbf{Matrix-based update form.}
We can rewrite Algorithm \ref{alg:stpp} in a compact matrix form:
\begin{equation}
\label{eq:matrix-update}
\begin{aligned}
\bX^{(t+1)} & = \bR\,\bigl(\bX^{(t)} - \gamma\,\bY^{(t)}\bigr),\\
\bY^{(t+1)} & = \bC\,\bY^{(t)} + \bG^{(t+1)} - \bG^{(t)},
\end{aligned}
\end{equation}
where $\bY^{(0)} = \bG^{(0)}$, and the matrices $\bR, \bC\in \reals^{n\times n}$ are as defined in \eqref{eq:rc}. The vectors $\bG^{(t)}$ collect the stochastic gradients across all agents at iteration~$t$. Such a representation highlights how the spanning trees govern the information flow and how the proposed method leverages both push and pull dynamics to achieve distributed optimization efficiently.

\section{Preliminary Results}
\label{sec:preliminary}
In this section, we present several preliminary results covering the introduction of useful analytical tools and convergence analysis, based on the properties of the spanning trees and their associated mixing matrices.

We begin by outlining the assumptions necessary for analyzing the convergence of the proposed STPP algorithm. These assumptions relate to the stochastic gradients, the smoothness of the local cost functions, and the strong convexity of the overall objective function.

To solve problem \eqref{eq:obj}, we assume each agent $i$ queries a stochastic oracle ($\cS\cO$) to obtain noisy gradient samples of the form $g_i(x_i,\xi_i)$, which satisfy the following conditions:
\begin{assumption}
    \label{a.var}
    For all $i \in \cN$ and $x\in\reals^p$, each random vector $\xi_i$ is independent, and  the stochastic gradients satisfy
\[
\begin{aligned}
   & \expect_{\xi_i\sim\cD_i} \brk{ g_i(x,\xi_i)| x}  = \nabla f_i(x),  \\
   & \expect_{\xi_i \sim \cD_i} \brk{ \norm{g_i(x,\xi_i) - \nabla f_i(x)}^2 |x}  \le \sigma^2,
\end{aligned}
\]
for some $\sigma^2 > 0$.
\end{assumption}

We make the following standard standing assumption on the individual objective functions $f_i$.
\begin{assumption}
    \label{a.smooth}
    The local cost functions $f_i(x): \reals^p \rightarrow \reals$ are lower bounded and have $L$-Lipschitz continuous gradients, i.e., for any $x,y \in \reals^p$,
        \[
        \norm{ \nabla f_i(x) - \nabla f_i(y)} \le L\norm{x - y}.
        \]
\end{assumption}

We also consider the following standard assumption regarding strong convexity.
\begin{assumption}
    \label{a.convex}
    The global objective function $f(x): \reals^p \rightarrow \reals$ is $\mu$-strongly convex, i.e., for any $x,y \in \reals^p$,
    \[
    f(y) \ge f(x) + \ip{\nabla f(x)}{y-x} + \frac{\mu}{2}\norm{y-x}^2.
    \]
\end{assumption}
Let $f^*:= \min_{x\in\reals^p}f(x)$ denote the minimum value of $f$ if Assumption \ref{a.smooth} holds. And let $x^*:=\arg\min_{x\in\reals^p} f(x)$ be the corresponding optimal solution if Assumption \ref{a.convex} holds.

\subsection{Supporting Lemmas}
In this subsection, we present several essential inequalities and lemmas that are instrumental in proving our main results.
\begin{lemma}
    \label{lem:sum_matrix}
    For an arbitrary set of $m$ matrices $\crk{\bA_i}_{i=1}^{m}$ with the same size, we have
    \[
    \norm{\sum_{i=1}^{m}\bA_i}_F^2 \le m \sum_{i=1}^{m} \norm{\bA_i}_F^2.
    \]
\end{lemma}
\begin{proof}
    See Lemma A.5 in \cite{you2024b} for reference. 
\end{proof}
\begin{lemma}
    \label{lem:sum_helper}
    For an arbitrary set of $p$ matrices $\crk{\bA_i}_{i=1}^{p}$ with the same size, and given $\alpha_1,\cdots, \alpha_p$ satisfying $\alpha_i > 0$ and $\sum_{i=1}^{p}\alpha_i \le 1$, we have
    \[
    \norm{\sum_{i=1}^{p}\bA_i}_F^2 \le \sum_{i=1}^{p} \frac{1}{\alpha_i} \norm{\bA_i}_F^2.
    \]
\end{lemma}
\begin{proof}
See Inequality (2) in \cite{song2024optimal} for reference.
\end{proof}
\begin{lemma}
    \label{lem:matrix_norm}
    Let $\bA$, $\bB$ be two matrices whose sizes match. Then,
    \[
    \norm{\bA\bB}_F \le \norm{\bA}_2\norm{\bB}_F.
    \]
\end{lemma}
\begin{proof}
    See Lemma A.6 in \cite{you2024b} for reference.
\end{proof}

\begin{lemma}
    \label{lem:martin}
    Let $A,B,C$ and $\alpha$ be positive constants and $T$ be a positive integer. Define 
    \[
    g(\gamma) = \frac{A}{\gamma (T+1)} + B\gamma + C\gamma^2.
    \]
    Then,
    \[
    \inf_{\gamma \in (0,\frac{1}{\alpha}]} g(\gamma) \le 2\prt{\frac{AB}{T+1}}^{\frac{1}{2}} + 2C^{\frac{1}{3}} \prt{\frac{A}{T+1}}^{\frac{2}{3}} + \frac{\alpha A}{T+1},
    \]
    where the upper bound could be achieved by choosing $\gamma = \min\crk{\prt{\frac{A}{B(T+1)}}^{\frac{1}{2}} , \prt{\frac{A}{C(T+1)}}^{\frac{1}{3}},\frac{1}{\alpha}}$
\end{lemma}
\begin{proof}
    See Lemma 26 in \cite{koloskova2021improved} for reference.
\end{proof}

\begin{lemma}
    \label{lem:sum_help}
    For non-negative series $\crk{a_t}_{t=0}^{T}$, $\crk{b_i}_{i=0}^{d-1}$ given integers $T, d$, we have
    \[
    \sum_{t=0}^{T} \sum_{m=\max\crk{0,t-d}}^{t-1} b_{t-m-1} a_m \le \sum_{i=0}^{d-1} b_i \sum_{t=0}^{T} a_t,
    \]
    where any illegal summation is set to be $0$, for example, $\sum_{m=\max\crk{0,t-d}}^{t-1}b_{t-m-1}  a_m = 0$ if $t=0$.
\end{lemma}

\begin{proof}
    We denote $i= t-m-1$. Then,
    \[
    \begin{aligned}
        & \sum_{t=0}^{T} \sum_{m=\max\crk{0,t-d}}^{t-1} b_{t-m-1} a_m  = \sum_{t=0}^{T} \sum_{i=0}^{\min\crk{t,d}-1} b_{i} a_{t-i-1}  \\
         = & \sum_{i=0}^{d-1} \sum_{t=i+1}^{T} b_i a_{t-i-1} = \sum_{i=0}^{d-1} b_i \sum_{t=0}^{T-i-1} a_t  \le \sum_{i=0}^{d-1} b_i \sum_{t=0}^{T} a_t.
    \end{aligned}
    \]
\end{proof}
\begin{lemma}
    \label{lem:independ_help}
    Consider three random variables $X$, $Y$, and $Z$. Assume that $X$ and $Z$ are independent, $Y$ and $Z$ are independent and also $Z$ and $(X,Y)$ are independent. Let $h$ and $g$ be functions such that the conditional expectation $\mathbb{E}[g(Y,Z) \mid Y] = 0$. Under these conditions, we have
    \[
    \expect \prt{h(X) g(Y,Z)} = 0.
    \]
\end{lemma}
\begin{proof}
    See Lemma A.8 in \cite{you2024b} for reference.
\end{proof}
\begin{lemma}
    \label{lem:gamma}
    Given a positive constant $d$, for any stepsize $\gamma \le \frac{1}{ 3 n d  L}$, the following inequality holds:
    \[
    \prt{1-\frac{\gamma n \mu }{4}}^{d} \ge \gamma n d  L.
    \]
    Specifically, in the subsequent lemmas, useful choices for the parameter $d$ are $d = \kappa d_\bC$ and $d = \kappa \prt{d_\bR + d_\bC}$.
    
\end{lemma}
\begin{proof}
    It holds by $\gamma \le \frac{1}{ 3 n d  L}$ that
    \begin{equation}
        \label{eq:lem:gamma-1}
        \begin{aligned}
            \frac{1}{d} & \ge  \gamma n \brk{ L + \frac{\mu}{4}+ \frac{\mu}{4d}} \ge \gamma n \brk{L + \frac{\mu}{4}  + \frac{\mu}{4d}- \frac{\gamma n  L \mu}{4}}.
        \end{aligned}
    \end{equation}
    Then, by rearranging terms, we have $-\frac{\gamma n \mu}{4} \ge \gamma n  L - \frac{\gamma^2 n^2  L \mu}{4} - \frac{1}{d} + \frac{\gamma n \mu}{4d}$ which implies
     $-\frac{\gamma n \mu}{4} \ge \frac{1}{d}\prt{\gamma n d L - 1}\prt{1 - \frac{\gamma n \mu}{4}}.$
    It follows from $0< 1 - \frac{\gamma n \mu}{4} < 1$ and the logarithm inequality $1-\frac{1}{x} \le \log(x) \le x - 1$ that
    \begin{equation}
        \label{eq:lem:gamma-4}
        \begin{aligned}
            & \log\prt{1-\frac{\gamma n \mu}{4}} \ge 1 -\frac{1}{1-\frac{\gamma n \mu}{4}} = \frac{- \frac{\gamma n \mu}{4}}{1-\frac{\gamma n \mu}{4}} \\
            \ge &  \frac{1}{d}\prt{\gamma n d  L - 1}\ge \frac{1}{d} \log\prt{\gamma n d  L }.
        \end{aligned}
    \end{equation}
    Taking exponents on both sides yields the desired result.
\end{proof}

\subsection{Convergence Analysis}
In this subsection, we present the key steps for analyzing the convergence of the STPP algorithm. We first introduce the matrices $\bPi_\bR:= \bI - \mone\bpi_\bR^\T$ and $\bPi_\bC:= \bI - \bpi_\bC\mone^\T$, and define another important term to capture the property of the spanning tree $\bA_k:= \bC^{k} - \bC^{k-1}$ for $k=1,\cdots ,d_{\bC}$ and $\bA_0 := \bI$. Intuitively speaking, $\prt{ \bpi_{\bR} - \mone}^\T \bA_{k}$ indicates the nodes with length $k$ to node 1 in the spanning tree $\cG_{\bC}$. Furthermore, we define
\[
\begin{aligned}
    \bar{\bX}^{(t)} & := \mone\bpi_\bR^\T \bX^{(t)}, \Delta \bar{\bX}^{(t)} := \bar{\bX}^{(t+1)} - \bar{\bX}^{(t)},\\
    \hat{\bX}^{(t)} &:= \bPi_\bR\bX^{(t)},\nabla \bar{\bF}^{(t)} := \nabla \bF^{(t)} - \nabla \bF(\bar{\bX}^{(t)}).
\end{aligned}
\]

Next, we express STPP as follows following Equation \eqref{eq:matrix-update}:
\begin{equation}
    \label{eq:matrixpp}
    \brk{ \begin{array}{c}
         \bX^{(t+1)}\\
        \bY^{(t+1)}  
    \end{array}} = 
    \brk{\begin{array}{cc}
         \bR & -\gamma \bR\\
        \mathbf{0} & \bC  
    \end{array}} \brk{\begin{array}{c}
         \bX^{(t)}\\
        \bY^{(t)}  
    \end{array}} + 
    \brk{\begin{array}{c}
         \mathbf{0}\\
        \bG^{(t+1)} -\bG^{(t)} 
    \end{array}}.
\end{equation}
By iteratively applying Equation \eqref{eq:matrixpp}, starting from time step $t$ and going back to the initial step $t=0$, we obtain the following expression:
\[
\begin{aligned}
    \brk{ \begin{array}{c}
         \bX^{(t)}\\
        \bY^{(t)}  
    \end{array}}
    & = \brk{ \begin{array}{cc}
         \bR & -\gamma \bR\\
        \mathbf{0} & \bC  
    \end{array}}^t\brk{ \begin{array}{c}
         \bX^{(0)}\\
        \bY^{(0)}  
    \end{array}} \\
    & + 
    \sum_{m=0}^{t-1}
    \brk{ \begin{array}{cc}
         \bR & -\gamma \bR\\
        \mathbf{0} & \bC  
    \end{array}}^{t-m-1}\brk{\begin{array}{c}
         \mathbf{0}\\
        \bG^{(m+1)} -\bG^{(m)} 
    \end{array}}.
\end{aligned}
\]
Note that for any integer $m>0$, we have
\[
\brk{ \begin{array}{cc}
         \bR & -\gamma \bR\\
        \mathbf{0} & \bC  
    \end{array}}^m = \brk{ \begin{array}{cc}
         \bR^m & -\gamma \sum_{j=1}^{m}\bR^j\bC^{m-j}\\
        \mathbf{0} & \bC^m  
    \end{array}}.
\]
Therefore, starting from the initial condition $\bX^{(0)} = \bG^{(0)}$, we can deduce the expressions for $\bX^{(t)}$ and $\bY^{(t)}$ as follows:
\begin{equation}
    \label{eq:pp1}
    \begin{aligned}
        \bX^{(t)} & = \bR^t\bX^{(0)}   -\gamma \sum_{j=1}^{t}\bR^j\bC^{t-j} \bG^{(0)} \\
        & - \gamma \sum_{m=0}^{t-2} \sum_{j=1}^{t-m-1}\bR^j\bC^{t-m-1-j} \brk{\bG^{(m+1)} - \bG^{(m)}},
    \end{aligned}
\end{equation}
\begin{equation}
    \label{eq:pp2}
    \bY^{(t)}  = \sum_{m=0}^{t-1} \bC^{t-m-1} \brk{\bG^{(m+1)} - \bG^{(m)}} + \bC^t \bG^{(0)}.
\end{equation}
Then, by multiplying $\bPi_\bC$ to both sides of Equation \eqref{eq:pp2} and invoking Lemma \ref{lem:Rk_d}, we have
\begin{equation}
    \label{eq:pp22}
    \begin{aligned}
    & \bPi_\bC\bY^{(t)} = \sum_{m=1}^{\min\{t,d_\bC\}}\bA_m \bG^{(t-m)} + \bPi_{\bC}\bG^{(t)}.\\
    \end{aligned}
\end{equation}

To support the convergence analysis, Lemma \ref{lem:var} provides an estimate for the variance of gradient estimator $\bG(\bX^{(t)},\bxi^{(t)})$. Following the result, Corollary \ref{lem:var0} and Corollary \ref{lem:var1} introduce useful tools.
\begin{lemma}
    \label{lem:var}
    Suppose Assumption \ref{a.var} holds. Then, for any index set $\cI \subseteq [n]$, we have
    \[
    \begin{aligned}
        & \expect\brk{\norm{\bfe_{\cI}^\T \bTh^{(t)}}_2^2 }  \le |\cI| \sigma^2. 
    \end{aligned}
    \]
\end{lemma}
\begin{proof}
    This result follows from Assumption \ref{a.var}.
\end{proof}

\begin{corollary}
    \label{lem:var0}
    Suppose Assumption \ref{a.var} holds. Then, for any given power $k\le d_\bC - 1$, we have 
    \[
    \expect\brk{\norm{\prt{\bC^{k}-\bpi_{\bC}\mone^\T}\bTh^{(t)}}_F^2} \le 2\prt{n-c_k}\sigma^2.
    \]
\end{corollary}
\begin{proof}
    Invoking Lemma \ref{lem:Rk_d}, we have 
    \[
    \begin{aligned}
        & \norm{\brk{\bC^{k}-\bpi_{\bC}\mone^\T}\bTh^{(t)}}_F^2 =   2\norm{\brk{\bfe_{\tilde{\cI}_{\bC,1}^{k}} - \mone}^\T\bTh^{(t)}  }_2^2 
    \end{aligned}
    \]
    Taking expectation on both sides, from Lemma \ref{lem:var} we have
    \[
    \expect\brk{\norm{\prt{\bC^{k}-\bpi_{\bC}\mone^\T}\bTh^{(t)}}_F^2 } \le 2\prt{n-c_k} \sigma^2.
    \]
\end{proof}
\begin{corollary}
    \label{lem:var1}
    Suppose Assumption \ref{a.var} holds. Then, we have
    \[
    \sum_{k=1}^{d_\bC}\expect \norm{\prt{\bpi_\bR - \mone}^\T\bA_{k}\bTh^{(t)}}_F^2  \le \prt{n-1 }\sigma^2.
    \]
\end{corollary}
\begin{proof}
     By Corollary \ref{lem:var0}, there holds
    \[
    \begin{aligned}
        & \sum_{k=1}^{d_\bC}\expect \norm{\brk{\bpi_\bR - \mone}^\T\bA_{k}\bTh^{(t)}}_F^2  =  \sum_{k=1}^{d_\bC}\expect \norm{\brk{\bfe_{\tilde{\cI}_{\bC,1}^{k}}^\T - \bfe_{\tilde{\cI}_{\bC,1}^{k-1}}^\T }\bTh^{(t)}}_F^2  \\
        & \le \sum_{k=1}^{d_\bC} \prt{c_{k} - c_{k-1} }\sigma^2 = \prt{c_{d_\bC} - c_0}\sigma^2 = \prt{n-1}\sigma^2.
    \end{aligned}
    \]
\end{proof}

The following lemma is crucial for obtaining tighter convergence results compared to prior works, specifically \cite{you2024b}.
\begin{lemma}
    \label{lem:Ak-barF}
    Suppose Assumption \ref{a.smooth} holds. We have
    \[
    \begin{aligned}
        \sum_{t=0}^{T}\expect \norm{\sum_{m=0}^{\min\crk{t,d_\bC}}\bpi_\bR^\T\bA_m\nabla \bar{\bF}^{(t-m)}}^2 & \le nL^2\sum_{t=0}^{T}\expect \norm{\bPi_\bR\bX^{(t)}}_F^2.
    \end{aligned}
    \]
    Furthermore, for $\gamma \le \frac{1}{3nd_\bC \kappa L}$, we have
    \[
    \begin{aligned}
        & \sum_{t=0}^{T}\prt{1-\frac{\gamma n \mu}{4}}^{T-t}\expect \norm{\sum_{m=0}^{\min\crk{t,d_\bC}}\bpi_\bR^\T\bA_m\nabla \bar{\bF}^{(t-m)}}^2 \\
        \le & n L^2 \sum_{t=0}^{T} \prt{1-\frac{\gamma n \mu}{4}}^{T-t-d_\bC} \norm{\bPi_\bR\bX^{(t)}}_F^2 \\
        & \prt{\le   \frac{L}{\gamma d_\bC  \kappa }\sum_{t=0}^{T} \prt{1-\frac{\gamma n \mu}{4}}^{T-t} \norm{\bPi_\bR\bX^{(t)}}_F^2}. 
    \end{aligned}
    \]
\end{lemma}
\begin{proof}
    See Appendix \ref{pf:lem:Ak-barF}.
\end{proof}
Lemma \ref{lem:xk-Ak-2} plays a crucial role for obtaining tighter convergence results under Assumption \ref{a.convex}.
\begin{lemma}
    \label{lem:xk-Ak-2}
    Suppose Assumption \ref{a.smooth} holds. Then, we have for $\gamma \le 1/\prt{3nd_\bC \kappa L}$ that
    \[
    \begin{aligned}
        & \sum_{t=0}^{T}\prt{1-\frac{\gamma n \mu}{4}}^{T-t} \norm{\bpi_\bR^\T \sum_{m=0}^{\min\crk{t,d_\bC}}\bA_m \nabla \bF(\bar{\bX}^{(t-m)})}^2 \\
        \le &   \frac{6 L}{\gamma  \kappa }\sum_{t=0}^{T} \prt{1-\frac{\gamma n \mu}{4}}^{T-t} \norm{\Delta \bar{\bX}^{(t)} }^2 + 6nc_\avg \prt{1-\frac{\gamma n \mu}{4}}^{T-d_\bC} \cdot \\
        & \quad \norm{ \nabla \bF^{(0)}}^2_F  + 3n^2\sum_{t=0}^{T} \prt{1-\frac{\gamma n \mu}{4}}^{T-t}\norm{\nabla f(x_1^{(t)})}^2.
    \end{aligned}
    \]
\end{lemma}
\begin{proof}
    See Appendix \ref{pf:lem:xk-Ak-2}.
\end{proof}

The following lemmas delineate the critical components for bounding the average expected squared norm of the gradient, i.e., $\frac{1}{T+1} \sum_{t=0}^{T}\expect\norm{\nabla f(x_1^{(t)})}^2$ under Assumptions \ref{a.var} and \ref{a.smooth}, and the optimality gap $\expect\norm{x^{(T)}_1 - x^*}^2$ under Assumptions \ref{a.var}-\ref{a.convex}.

First, Lemma \ref{lem:X_diff} introduces the primary bounds on the expressions $\sum_{t=0}^{T}\expect \norm{\Delta\bar{\bX}^{(t)} }_F^2$ and $\sum_{t=0}^{T}\prt{1-\frac{\gamma n \mu}{4}}^{T-t}\expect \norm{\Delta\bar{\bX}^{(t)} }_F^2$.
\begin{lemma}
    \label{lem:X_diff}
    Suppose Assumption \ref{a.var} and \ref{a.smooth} hold. For $\gamma \le 1/(6 n\sqrt{d_\bC c_\avg} L) $, we have the following inequality:
    \[
    \begin{aligned}
        &\sum_{t=0}^{T} \expect\norm{\Delta\bar{\bX}^{(t)} }_F^2  \le  10\gamma^2 n^2 \sigma^2\prt{T+1} + 10 \gamma^2 n^2 L^2 \sum_{t=0}^{T} \expect \norm{\hat{\bX}^{(t)}}_F^2\\
        &\qquad + 20 \gamma^2 n^2 c_{\avg}\norm{\nabla \bF^{(0)}}_F^2 + 10\gamma^2 n^3 \sum_{t=0}^{T} \expect\norm{\nabla f(x_1^{(t)})}^2.
    \end{aligned}
    \]
    In addition, for $\gamma \le \frac{1}{100 n d_\bC\kappa L}$ and $T > 2d_\bC$, we have
    \[
    \begin{aligned}
        & \sum_{t=0}^{T} \prt{1-\frac{\gamma n \mu}{4}}^{T-t} \expect \frac{1}{n}\norm{ \Delta\bar{\bX}^{(t)} }^2_F \\
        &\le    \frac{24\gamma \sigma^2}{ \mu} + 36\gamma^2 nc_\avg \prt{1-\frac{\gamma n \mu}{4}}^{T-d_\bC} \norm{ \nabla \bF^{(0)}}^2_F \\
        & + \sum_{t=0}^{T} \prt{1-\frac{\gamma n \mu}{4}}^{T-t} \expect\brk{ \frac{6\gamma L}{  d_\bC\kappa } \norm{\hat{\bX}^{(t)}}_F^2 + 18\gamma^2 n^2 \norm{\nabla f(x_1^{(t)})}^2 } .
    \end{aligned}
    \]
\end{lemma}
\begin{proof}
    See Appendix \ref{pf:lem:X_diff}.
\end{proof}

Lemma \ref{lem:PiX} provides the bounds for $\sum_{t=0}^{T}\expect \norm{\hat{\bX}^{(t)}}_F^2$ and $\sum_{t=0}^{T}\prt{1-\frac{\gamma n \mu}{4}}^{T-t}\expect \norm{\hat{\bX}^{(t)}}_F^2$, respectively.
\begin{lemma}
    \label{lem:PiX}
    Suppose Assumption \ref{a.var} and \ref{a.smooth} hold. For $\gamma \le 1/\prt{100 n\sqrt{d_\bR d_\bC r_\avg c_\avg} L } $, we have the following inequality:
    \[
    \begin{aligned}
        & \sum_{t=0}^{T}\expect\norm{\hat{\bX}^{(t)}}_F^2  \le 40 \gamma^2 n^3 d_\bR r_\avg \sum_{t=0}^{T} \expect\norm{\nabla f(x_1^{(t)})}^2 \\
        & \quad+  40 \gamma^2 n^2 \min\crk{d_\bR, d_\bC} r_\avg c_\avg \brk{\sigma^2\prt{T+1} + 2\norm{\nabla \bF^{(0)}}_F^2}.\\
    \end{aligned}
    \]
    Furthermore, for $\gamma \le 1/\prt{100 n\max\crk{d_\bR,d_\bC}\kappa L}$, we have
    \[
    \begin{aligned}
        & \sum_{t=0}^{T}\prt{1-\frac{\gamma n \mu}{4}}^{T-t} \expect \norm{\hat{\bX}^{(t)}}_F^2 \le  \frac{96 \gamma n \min\crk{d_\bC,d_\bC} r_\avg c_\avg \sigma^2}{ \mu} \\
        & + 34 \gamma^2 n^3 d_\bR \sum_{t=0}^{T} \prt{1-\frac{\gamma n \mu}{4}}^{T-t - d_\bC}\norm{\nabla f(x_1^{(t)})}^2 \\
        & + 104 \gamma^2 n^2 \prt{\min\crk{d_\bR,d_\bC}}^2\prt{d_\bR + d_\bC} \prt{1-\frac{\gamma n \mu}{4}}^{T-d_\bR -d_\bC}\cdot\\
        & \qquad \norm{ \nabla \bF^{(0)}}_F^2 + 2000 \gamma^2 n^2 d_\bR  \sigma^2 .
    \end{aligned}
    \]
\end{lemma}
\begin{proof}
See Appendix \ref{pf:lem:PiX}.
    
\end{proof}

The design of STPP implies that there is an inherent delay in the transmission of information from the nodes far away from the root node to the root. As the information traverses through the trees $\cG_{\bR}$ and $\cG_{\bC}$, the delay becomes evident. Specifically, for nodes with distance $k$ to node $1$ on $\cG_{\bC}$, their information requires an additional $k$ iterations to successfully reach and impact the node 1, as demonstrated in Lemma \ref{lem:independ1} below.
\begin{lemma}
    \label{lem:independ1}
    Suppose Assumption \ref{a.var} holds. Then, for any given layer $i$ in the spanning tree $\cG_{\bC}$, i.e., $i \in \crk{0,1,2, \cdots, d_\bC}$ and time step $t > i $, we have
    \[
    \begin{aligned}
        & \expect\ip{ \nabla f(x_1^{(t)})}{ \bfe_{\cI_{\bC,1}^{i}}^\T\bTh^{(t-i)}} = 0, \  \expect\ip{ x_1^{(t)} - x^*}{ \bfe_{\cI_{\bC,1}^{i}}^\T\bTh^{(t-i)}} = 0.
    \end{aligned}
    \]
\end{lemma}

\begin{proof}
See Appendix \ref{pf:lem:independ1}.
\end{proof}

Building on the preceding lemmas, we are in a position to establish the key auxiliary convergence results for the STPP algorithm, which involves bounding the average expected squared norm of the gradient, as detailed in Lemma \ref{lem:smooth}. 
\begin{lemma}
    \label{lem:smooth}
    Suppose Assumption \ref{a.var} and \ref{a.smooth} hold. For $\gamma \le  1/\prt{100n\sqrt{d_\bR d_\bC r_\avg c_\avg} L}$, we have
    \[
    \begin{aligned}
        &\frac{1}{T+1} \sum_{t=0}^{T}\expect \norm{\nabla f(x_1^{(t)})}^2 \le \frac{8\Delta_f}{\gamma n \prt{T+1}} + 40 \gamma L \sigma^2 \\
        & \quad + 2400 nd_\bC r_\avg c_\avg \gamma^2  L^2 \sigma^2 +  \frac{48 c_\avg}{n\prt{T+1}}\norm{\nabla \bF^{(0)}}^2_F.
    \end{aligned}
    \]
\end{lemma}

\begin{proof}
See Appendix \ref{pf:lem:smooth}.
\end{proof}
 Under  strong convexity, we bound the optimality gap.
\begin{lemma}
    \label{lem:convex}
    Suppose Assumptions \ref{a.var}, \ref{a.smooth} and \ref{a.convex} hold. For $\gamma \le 1/\prt{1000\max\crk{d_\bR,d_\bC} r_\avg c_\avg \kappa L}$, we have 
    \[
    \begin{aligned}
        & \expect \norm{x_1^{(T)} - x^*}^2 \le \prt{1 - \frac{\gamma n \mu}{4}}^{T} \brk{\norm{x_1^{(0)} - x^*}^2 + \frac{1}{n L^2} \norm{\nabla \bF^{(0)}}_F^2 } \\
        & \quad +  \frac{2000\gamma \kappa \sigma^2}{ \mu}  + 40000 \gamma^2 n \min\crk{d_\bR,d_\bC} r_\avg c_\avg \kappa^2 \sigma^2.
    \end{aligned}
    \]
\end{lemma}

\begin{proof}
See Appendix \ref{pf:lem:convex}. 
\end{proof}

\section{Main Results}
\label{sec:convergence}

The main convergence properties of STPP are summarized in the following theorems. The first result applies to the general non-convex case, while the second one assumes strong convexity of $f$.

\begin{theorem}
    \label{thm:smooth}
    Consider the STPP algorithm (Algorithm \ref{alg:stpp}) implemented on the spanning tree graphs $\cG_{\bR}$ and $\cG_{\bC}$. Assume Assumption \ref{a.var} and Assumption \ref{a.smooth} hold. Given the number of iterations $T$, define the stepsize $\gamma$ as
    \[
    \begin{aligned}
        & \gamma = \min\left\{\frac{1}{100n\sqrt{d_\bR d_\bC r_\avg c_\avg} L} , \prt{\frac{\Delta_f}{5nL\sigma^2(T+1)}}^{\frac{1}{2}}, \right.\\
        &\qquad \left. \prt{\frac{\Delta_f}{300n^2d_\bC r_\avg c_\avg L\sigma^2(T+1)}}^{\frac{1}{3}}  \right\}. 
    \end{aligned}
    \]  
    Then, the following convergence result holds:
    \begin{equation}
        \label{eq:thm:smooth}
        \begin{aligned}
            &\frac{1}{T+1} \sum_{t=0}^{T}\expect \norm{\nabla f(x_1^{(t)})}^2 \le 100\prt{\frac{\Delta_f L \sigma^2}{ n \prt{T+1}}}^{\frac{1}{2}} \\
            &  + 300 \prt{n \max\crk{d_\bR, d_\bC} r_\avg c_\avg L^2 \sigma^2}^{\frac{1}{3}} \prt{\frac{\Delta_f}{n(T+1)}}^{\frac{2}{3}} \\
            &  + \frac{800\sqrt{d_\bR d_\bC r_\avg c_\avg}\Delta_f}{T+1} +  \frac{50 c_\avg}{n\prt{T+1}}\norm{\nabla \bF^{(0)}}^2_F,
        \end{aligned}
    \end{equation}
    where $\Delta_f := f(x_1^{(0)}) - f^*$, $d_\bR, d_\bC$ represent the  diameters of the graphs $\cG_\bR$ and $\cG_\bC$, respectively, and $r_\avg, c_\avg $ represent the average  distances of each node to the root node 1.
\end{theorem}
 \begin{proof}
    Invoking Lemmas \ref{lem:martin} and \ref{lem:smooth}, where 
    \begin{equation}
        \begin{aligned}
        &\alpha = 100n\sqrt{d_\bR d_\bC r_\avg c_\avg} L, A= \frac{8\Delta_f}{n}, \\
        & B = 40 L\sigma^2,  C = 2400 nd_\bC r_\avg c_\avg L^2 \sigma^2,
        \end{aligned}
    \end{equation}
    we get the desired result.
\end{proof}
\begin{remark}
    Based on the convergence rate in Equation \eqref{eq:thm:smooth}, we derive that when $T = \cO\prt{n\prt{\max\crk{d_\bR, d_\bC} r_\avg c_\avg}^{2}}$, the term $\cO(\frac{1}{\sqrt{nT}})$ dominates the remaining terms up to a constant scalar. This implies that STPP achieves linear speedup after $\cO\prt{n\prt{\max\crk{d_\bR, d_\bC} r_\avg c_\avg}^{2}}$ transient iterations.
\end{remark}

\begin{remark}
    By the definitions of $d_\bR, d_\bC, r_\avg$ and $c_\avg$,  they are all bounded above by $n$. Thus we can derive a polynomial upper bound for the transient iterations, specifically $\cO(n^7)$. This upper bound is significantly  better than the worst-case performance of existing decentralized  stochastic gradient algorithms, such as those in \cite{assran2019stochastic, liang2023towards}. 
\end{remark}

\begin{theorem}
    \label{thm:convex}
    Consider the STPP algorithm (Algorithm \ref{alg:stpp}) implemented on the spanning tree graphs $\cG_{\bR}$ and $\cG_{\bC}$. Assume Assumption \ref{a.var}, \ref{a.smooth} and \ref{a.convex} hold. Given the number of iterations $T$, define the stepsize $\gamma$ as
    \[
    \begin{aligned}
        & \gamma = \min \left\{ \frac{16\log\prt{n(T+1)}}{n\mu(T+1)}, \frac{1}{1000n\max\crk{d_\bR, d_\bC} r_\avg c_\avg \kappa L} \right\}.
    \end{aligned}
    \]
    Then, the following convergence result holds:
    \begin{equation}
        \label{eq:thm:convex}
        \begin{aligned}
            & \expect\norm{x_1^{(T)} - x^*}^2 \le \frac{32000\sigma^2 \log\prt{n(T+1)} }{n(T+1)\mu^2} \\
            & \quad + \frac{10240000\max\crk{d_\bR,d_\bC} r_\avg c_\avg \kappa^2 \sigma^2 }{n(T+1)^2\mu^2} \\
            & \quad + \max \left\{ \exp\prt{-\frac{T}{4000\max\crk{d_\bR,d_\bC} r_\avg c_\avg \kappa^2}}, \right. \\
            &   \quad \left. \frac{40}{n(T+1)^2}\right\}\prt{\norm{x_1^{(0)} - x^*}^2 + \frac{1}{n L^2} \norm{\nabla \bF^{(0)}}_F^2 },
        \end{aligned}
    \end{equation}
    where $x^*:= \arg\min_{x} f(x)$ and $\kappa:=\frac{L}{\mu}$ is the conditional number.
\end{theorem}
\begin{proof}
    Given the choice of $\gamma$,
    we get
    \begin{equation}
        \label{eq:thm:convex-2}
        \begin{aligned}
            & \prt{1 - \frac{\gamma n \mu}{4}}^{T} \le \max\bigg\{\prt{1-\frac{1}{4000\max\crk{d_\bR,d_\bC} r_\avg c_\avg \kappa^2}}^{T},\\
            & \qquad \qquad \prt{1-\frac{4\log\prt{n(T+1)}}{T+1}}^{T}\bigg\} \\
            & \le \max\crk{\exp\prt{-\frac{T}{4000\max\crk{d_\bR,d_\bC} r_\avg c_\avg \kappa^2}}, \frac{40}{n(T+1)^2}},
        \end{aligned}
    \end{equation}
    where we use the fact $(1-\frac{1}{x})^{x} \le e^{-1}$ and $1-x \le \exp(-x)$ for any $x \in \reals_{+}$ in the last inequality.
    Then, invoking Lemma \ref{lem:convex} yields the desired result.
\end{proof}
\begin{remark}
    The convergence rate in Equation \eqref{eq:thm:convex} implies that 
    \[
    \begin{aligned}
        &\expect \norm{x_1^{(T)} - x^*}^2 \le \tilde{\cO} \left\{\frac{1}{nT} + \frac{\max\crk{d_\bR,d_\bC} r_\avg c_\avg }{nT^2} \right. \\
        & \left.  + \exp\prt{-\frac{T}{\max\crk{d_\bR,d_\bC} r_\avg c_\avg }} \right\},
    \end{aligned}
    \]
    where $\tilde{\cO}$ hides the constants and polylogarithmic factors.
    The transient time can be estimated as $\tilde{\cO}(\max\crk{d_\bR,d_\bC} r_\avg c_\avg )$. This corresponds to the number of iterations required before the term $\cO(\frac{1}{nT})$ dominates the  remaining terms.
\end{remark}
\begin{remark}
    Similar to the nonconvex case, since $d_\bR$, $d_\bC$, $r_\avg$ and $c_\avg$ are all bounded above by $n$, we can derive a polynomial upper bound for the transient iterations, specifically $\tilde{\cO}(n^3)$. This upper bound is significantly better than the worst-case performance of existing  algorithms, such as those in \cite{assran2019stochastic, liang2023towards}. 
\end{remark}

Compared to prior works such as \cite{alghunaim2022unified, koloskova2021improved, pu2021sharp}, which primarily rely on the spectral gap of the  mixing matrix to characterize the transient time, we adopt a different approach by using the diameter and average  distances of the spanning trees. 
Different from the general  observation that ``sparser graphs lead to slower convergence'', our algorithm and analysis offer a more refined conclusion: the convergence rate is explicitly determined by the diameter and average  distances of the spanning trees. Specifically, when the spanning trees drawn from the graph exhibit smaller diameters and average distances, the transient time is reduced, leading to faster convergence (see Table \ref{tab:smooth} and Table \ref{tab:convex}).

\section{Numerical Experiments}
\label{sec:exp}

In this section, we validate the previous theoretical results via numerical experiments. First, we examine the performance of the STPP algorithm and compare it with other popular  distributed stochastic gradient algorithms under several commonly  considered communication graphs 
on a standard nonconvex logistic regression  problem. Then, we apply the STPP algorithm on  deep learning tasks to show  its good performance in both training loss and test accuracy with respect to the number of iterations. The codes used to generate the figures in this section are available on github\footnote{\href{https://github.com/ryou98/SpanningTreePushPull}{https://github.com/ryou98/SpanningTreePushPull}}.

\textbf{Logistic Regression.} We compare the performance of STPP against other algorithms listed in Table \ref{tab:smooth} for logistic regression with non-convex regularization \cite{song2022communication,you2024b}. The objective functions $f_i: \mathbb{R}^p\rightarrow \mathbb{R}$ are given by
\[
f_i(x) := \frac{1}{J} \sum_{j = 1}^{J} \ln\prt{1+ \exp(-y_{i,j}h_{i,j}^\T x)} + R\sum_{k = 1}^{p} \frac{x_{[k]}^2}{1+ x_{[k]}^2},
\]
where $x_{[k]}$ is the $k$-th element of $x$, and $\crk{(h_{i,j}, y_{i,j})}_{j=1}^J$ represent the local data kept by node $i$. 
To control the data heterogeneity across the nodes, we first let each node $i$ be associated with a local logistic regression model with parameter $\tilde{x}_i$ generated by $\tilde{x}_i = \tilde{x} + v_i$, where $\tilde{x} \sim \cN(0,\bI_p)$ is a common random vector, and $v_i\sim\cN(0, \sigma_h^2\bI_p)$ are random vectors generated independently. Therefore, $\{v_i\}$ decide the dissimilarities between $\tilde{x}_i$, and larger $\sigma_h$ generally amplifies the difference. After fixing $\crk{\tilde{x}_i}$, local data samples are generated that follow distinct distributions. For node $i$, the feature vectors are generated as $h_{i,j} \sim \cN(0, \bI_p)$, and $z_{i,j}\sim\cU(0,1)$. Then, the labels $y_{i,j}\in \crk{-1,1}$ are set to satisfy $z_{i,j}\le 1 + \exp(-y_{i,j}h_{i,j}^\T \tilde{x}_i)$. In the simulations, the parameters are set as follows: $n=20$, $p=400$, $J = 500$, $R = 0.01$, and $\sigma_h = 0.2$. All the algorithms initialize with the same stepsize $\gamma = 0.4$, except STPP, which employs a modified stepsize $\gamma/n$. Such an adjustment is due to STPP's update mechanism, which incorporates a tracking estimator that effectively accumulates $n$ times the averaged stochastic gradients as the number of iterations increases. This can also be seen from the stepsize choice in Theorem \ref{thm:smooth}.\footnote{Note that this particular configuration results in slower convergence for STPP during the initial iterations, which can be improved by using larger initial stepsizes.} Additionally, we implement a stepsize decay of 80\% every $300$ iterations to facilitate convergence. 

\begin{figure}[!ht]
    \begin{center}
    \centerline{\includegraphics[width=0.8\columnwidth]{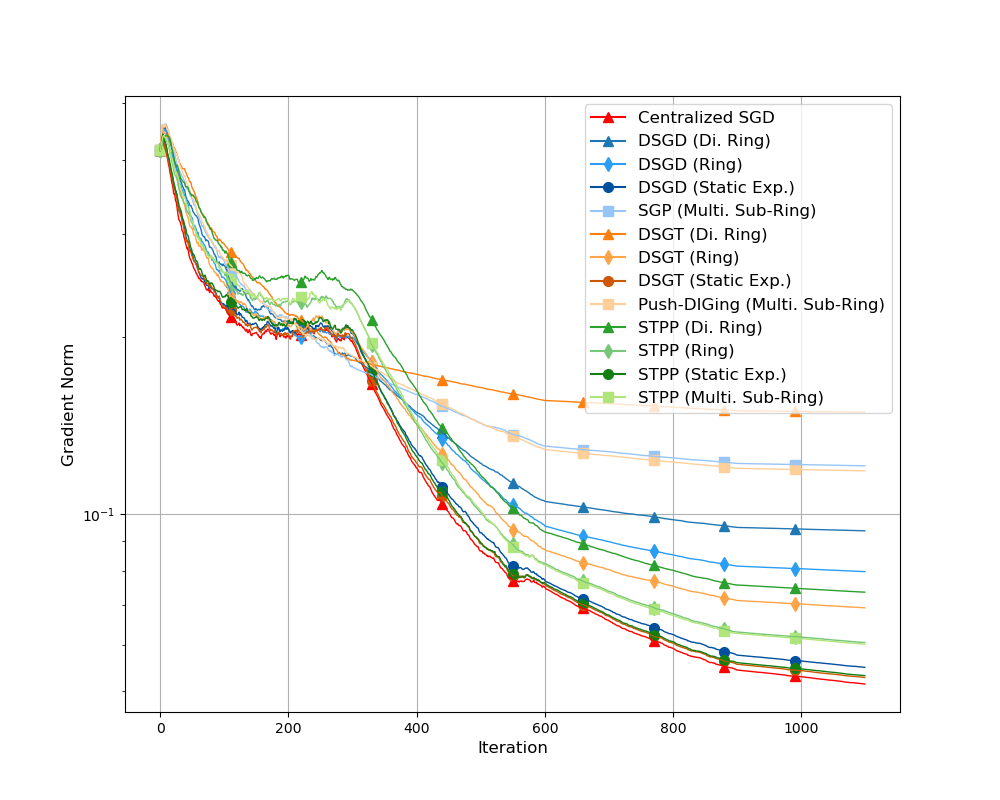}}
    \caption{Performance of different algorithms using various communication graphs for logistic regression with nonconvex regularization.}
    \label{fig:logistic}
    \end{center}
\end{figure}

Figure \ref{fig:logistic} illustrates the performance of the algorithms under various communication graphs. It is observed that STPP outperforms the other algorithms in sparse networks (e.g., Di. Ring and Ring) and non-regular  networks (e.g., Multi. Sub-Ring). For  dense networks (e.g., Static Exp.), STPP still achieves competitive performance.  

\begin{figure}[htbp]
	\centering
	\subfloat{\includegraphics[width=0.8\linewidth]{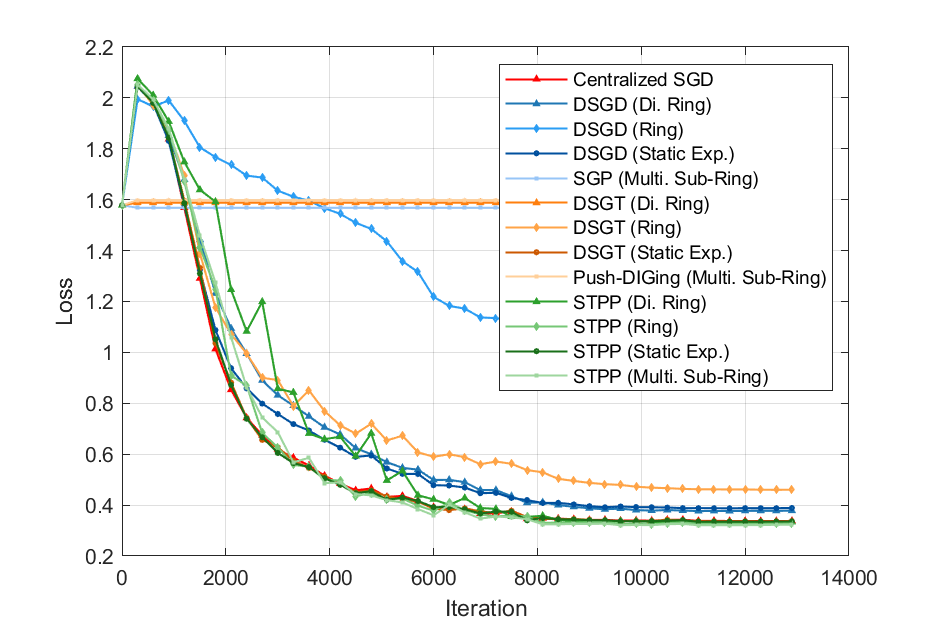}}
    \hfill
	\subfloat{\includegraphics[width=0.8\linewidth]{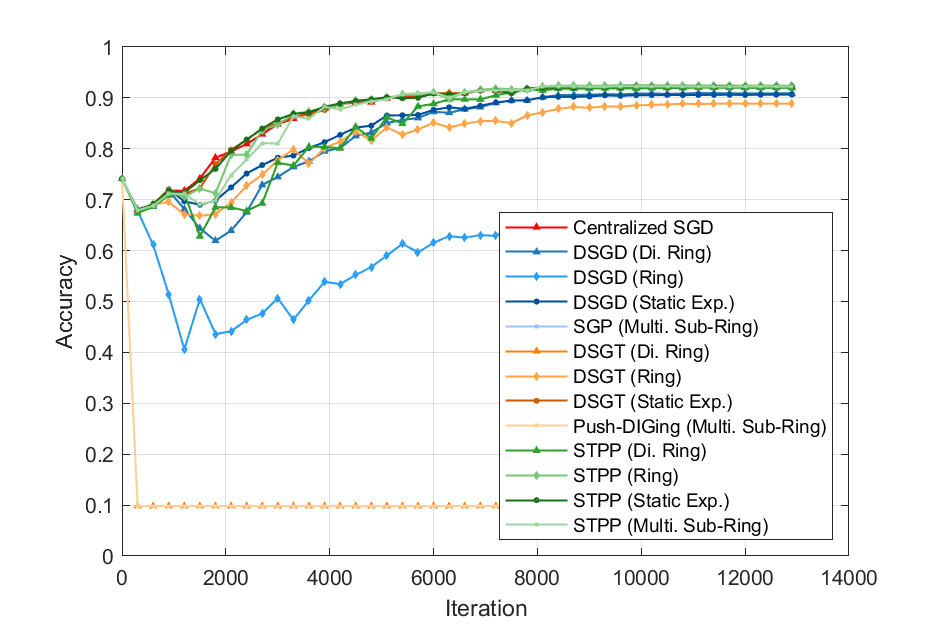}}
	\caption{\label{fig:mnist} Training loss and test accuracy of different algorithms for training CNN on MNIST.  For the training loss, we plot the curves of the divergentalgorithms (DSGT (Di. Ring), SGP (Multi. Sub-Ring) and Push-DIGing (Multi. Sub-Ring) with their first non-NaN values in the records.)}
	\label{fig:network}
\end{figure}

\textbf{Deep Learning.} We apply STPP and other algorithms to solve the image classification task with CNN over MNIST dataset \cite{lecun2010mnist}. We run all experiments on a server with eight Nvidia RTX 3090 GPUs. The network contains two convolutional layers with max pooling and ReLu and two feed-forward layers. In particular, we consider a heterogeneous data setting, where data samples are sorted based on their labels and partitioned among the agents. 
The local batch size is set to $8$ with $24$ agents in total. The learning rate is $0.01$ for all the algorithms except STPP (which employs a modified stepsize $\gamma/n$) for fairness. Additionally, the starting model is enhanced by pre-training using the SGD optimizer on the entire MNIST dataset for several iterations. The results are obtained by averaging over 3 independent random experiments. Figure \ref{fig:mnist} shows the  training loss and test accuracy curves. The results indicate that STPP outperforms the other algorithms in sparse networks (e.g., Di. Ring and Ring) and non-regular  networks (e.g. Multi. Sub-Ring). For dense networks (e.g., Static Exp.), STPP still achieves competitive performance.

\section{Conclusions}
\label{sec:conclusion}  

This paper introduces a novel distributed learning algorithm, \emph{Spanning Tree Push-Pull} (STPP), and provide comprehensive convergence analysis for both smooth nonconvex and strongly convex objectives. By leveraging two spanning trees derived from the underlying communication graph, STPP substantially improves upon existing decentralized optimization methods. In particular, the transient iteration complexity depends explicitly on the diameter and average distance of the spanning trees, thereby offering deeper insight into the relationship between network structure and convergence performance. Numerical experiments corroborate the theoretical findings.

\appendix

\section{Proofs}

\subsection{Proof of Lemma \ref{lem:Rk}}
\label{pf:lem:Rk}
 We prove the lemma by induction. First, it is obvious that $ \bR = \bZ_1$ by the definition of $\bR$:
    \[
    \begin{aligned}
        & \brk{\bR}_{i,j} = 1 \text{ iff } i \in \cI_{\bR,j}^{1} \text{ or } i= j =1\\
        \text{iff } & \brk{\bfe_{\cI_{\bR,j}^{1} }}_i = 1 \text{ or } i= j =1 \text{ iff } \brk{\bZ}_{i,j} = 1.
    \end{aligned}
    \]
    
    Assume the statement holds for $k = m$. Then, for $k = m+1$, we have
    \[
    \bR^{m+1} = \bR^{m}\bR = \bZ_{m}\bZ_{1}.
    \]
    Denote $\brk{\bZ_{m}\bZ_{1}}_i$ as the $i$-th column of matrix $\bZ_{m}\bZ_{1}$. To establish the result, we only need to demonstrate that the two matrices, $\bR^{m+1}$ and $\bZ_{m+1}$, share the same column entries. For the first column,
    \[
    \brk{\bZ_{m}\bZ_{1}}_1 = \bZ_{m}\brk{\bZ}_1 = \bZ_{m}\bfe_{\cI_{\bR,1}^{1}} =\bfe_{\tilde{\cI}_{\bR,1}^{m}} +  \sum_{i\in \cI_{\bR,1}^{1}} \bfe_{\cI_{\bR,i}^{m}}. 
    \]
    For the $p$-th element in $\brk{\bZ_{m}\bZ_{1}}_1$, it equals $1$ if and only if $p \in \tilde{\cI}_{\bR,1}^{m}$, or there exists a node $j \in \cI_{\bR,1}^1$ such that $p\in \cI_{\bR,j}^{m}$, or equivalently, the distance between node $p$ and node $1$ is no larger than $m+1$. Hence, we have $\brk{\bZ_{m}\bZ_{1}}_1 = \bfe_{\tilde{\cI}_{\bR,1}^{m+1}}$.
    
    For the $j$-th column where $j>1$, there holds
    \[
    \brk{\bZ_{m}\bZ_{1}}_j = \bZ_{m}\bfe_{\cI_{\bR,j}^{1}} = \sum_{i\in \cI_{\bR,j}^1} \bfe_{\cI_{\cR,i}^{m}}.
    \]
    Similarly, the $p$-th element in the $j$-th column $\brk{\bZ_{m}\bZ_{1}}_j$ equals $1$ if and only if there exists a node $i\in \cI_{\bR,j}^1$, such that  $p\in \cI_{\bR,j}^{m}$. Equivalently, this indicates that  the distance between node $p$ and node $j$ is exactly $m+1$. Thus, we have $\brk{\bZ_{m}\bZ_{1}}_j = \bfe_{\cI_{\bR,j}^{m+1}}$ which implies that $\bR^{m+1} = \bZ_{m+1}$, and we conclude that $\bR^{k} = \bZ_{k}$ for any positive index $k$.

\subsection{Proof of Lemma \ref{lem:Ak-barF}}
\label{pf:lem:Ak-barF}

From Lemma \ref{lem:Rk}, we have, by denoting $\bU_{t} := \norm{\sum_{m=0}^{\min\crk{t,d_\bC}}\bpi_\bR^\T\bA_m\nabla \bar{\bF}^{(t-m)}}^2$,
    \begin{equation}
        \label{eq:lem:Ak-barF-1}
        \begin{aligned}
            & \bU_{t} =   \langle\sum_{m=0}^{\min\crk{t,d_\bC}}\bpi_\bR^\T\bA_m\nabla \bar{\bF}^{(t-m)},\sum_{m=0}^{\min\crk{t,d_\bC}}\bpi_\bR^\T\bA_m\nabla \bar{\bF}^{(t-m)}\rangle \\
            = & \sum_{p=0}^{\min\crk{t,d_\bC}}\sum_{q=0}^{\min\crk{t,d_\bC}} \sum_{i\in \cI_{\bC,1}^{p}} \sum_{j\in \cI_{\bC,1}^{q}} \cdot \\
            & \ip{\nabla f_{i}(x_i^{(t-p)}) - \nabla f_{i}(x_1^{(t-p)})}{\nabla f_{j}(x_j^{(t-q)}) - \nabla f_{j}(x_1^{(t-q)})}.
        \end{aligned}
    \end{equation}
    Then, invoking Assumption \ref{a.smooth}, we have
    \begin{equation}
        \label{eq:lem:Ak-barF-2}
        \begin{aligned}
            & \sum_{t=0}^{T}\expect \bU_t \le  \sum_{t=0}^{T} \sum_{p=0}^{\min\crk{t,d_\bC}}\sum_{q=0}^{\min\crk{t,d_\bC}} \sum_{i\in \cI_{\bC,1}^{p}} \sum_{j\in \cI_{\bC,1}^{q}} \frac{L^2}{2}\cdot \\
            & \quad \brk{\expect  \norm{x_i^{(t-p)} - x_1^{(t-p)}}_2^2 + \expect\norm{x_j^{(t-q)} - x_1^{(t-q)} }_2^2} \\
            = &  \sum_{p=0}^{\min\crk{T,d_\bC}}\sum_{q=0}^{\min\crk{T,d_\bC}} \sum_{i\in \cI_{\bC,1}^{p}} \sum_{j\in \cI_{\bC,1}^{q}} \frac{L^2}{2} \cdot \\
            & \quad \sum_{t=p}^{T}  \brk{\expect  \norm{x_i^{(t-p)} - x_1^{(t-p)}}_2^2 + \expect\norm{x_j^{(t-q)} - x_1^{(t-q)} }_2^2}.
        \end{aligned}
    \end{equation}
    Noticing that $\sum_{p=0}^{\min\crk{T,d_\bC}} \sum_{i\in \cI_{\bC,1}^{p}} 1 \le n$, we get
    \begin{equation}
        \label{eq:lem:Ak-barF-3}
        \begin{aligned}
            & \sum_{t=0}^{T}\expect \bU_t
            \le  n L^2\sum_{p=0}^{\min\crk{T,d_\bC}} \sum_{i\in \cI_{\bC,1}^{p}}\cdot \\
            & \qquad \sum_{t=0}^{T}\expect \norm{x_i^{(t)} - x_1^{(t)}}_2^2 \le  n L^2 \sum_{t=0}^{T}\expect\norm{\bPi_\bR\bX^{(t)}}_F^2.
        \end{aligned}
    \end{equation}
    Furthermore, invoking Lemma \ref{lem:sum_helper} where $\alpha_m = \prt{1-\frac{\gamma n \mu}{4}}^{d_\bC - m}$, we get
    \begin{equation}
        \label{eq:lem:xk-Ak-1-2}
        \begin{aligned}
            \hat{\bU}:=& \sum_{t=0}^{T}\prt{1-\frac{\gamma n \mu}{4}}^{T-t} \norm{\sum_{m=0}^{\min\crk{t,d_\bC}}\bpi_\bR^\T\bA_m\nabla \bar{\bF}^{(t-m)}}^2 \\
            \le &   \sum_{t=0}^{T} \sum_{m=0}^{\min\crk{t,d_\bC}}  \prt{1-\frac{\gamma n \mu}{4}}^{k_{t,m}} \norm{\bpi_\bR^\T\bA_m}^2 L^2 \norm{\hat {\bX}^{(t-m)}}_F^2 \\
            =  &  \sum_{t=0}^{T} \sum_{j=t - \min\crk{t,d_\bC}}^{t}  \prt{1-\frac{\gamma n \mu}{4}}^{h_j} \norm{\bpi_\bR^\T\bA_{t-j}}^2 L^2 \norm{\hat{\bX}^{(j)}}_F^2 ,
        \end{aligned}
    \end{equation}
    where $k_{t,m} := T-t + m - d_\bC$ and $h_{j} := T-j - d_\bC$. Then,
    \[
    \begin{aligned}
        \hat{\bU}\le & \sum_{j=0}^{T} \sum_{t=j}^{j+d_\bC}  \prt{1-\frac{\gamma n \mu}{4}}^{T-j - d_\bC} \norm{\bpi_\bR^\T\bA_{t-j}}^2 L^2 \norm{\hat{ \bX}^{(j)}}_F^2 \\
            \le & n L^2 \sum_{t=0}^{T} \prt{1-\frac{\gamma n \mu}{4}}^{T-t-d_\bC} \norm{\hat{\bX}^{(t)}}_F^2.
    \end{aligned}
    \]
    We obtain the desired result by noting Lemma \ref{lem:gamma}.

\subsection{Proof of Lemma \ref{lem:xk-Ak-2}}
\label{pf:lem:xk-Ak-2}
Notice that, by Equation \eqref{eq:pp22},
    \begin{equation}
    \label{eq:lem:xk-Ak-2-1}
        \begin{aligned}
            & \bpi_\bR^\T \sum_{m=0}^{\min\crk{t,d_\bC}}\bA_m \nabla \bF(\bar{\bX}^{(t-m)}) \\
            = & -\sum_{m=0}^{\min\crk{t,d_\bC} - 1}\sum_{i\notin \tilde{\cI}_{\bC,1}^{m}} \prt{\nabla f_i(x_1^{(t-m)}) - \nabla f_i(x_1^{(t-m-1)})} \\
            & + \bpi_\bR^\T \bPi_\bC\bC^t \nabla \bF(\bX^{(0)})+ \bpi_\bR^\T\bpi_\bC \mone^\T \nabla \bF(\bar{\bX}^{(t)}).\\
        \end{aligned}
    \end{equation}
    It holds that $\sum_{i\notin \tilde{\cI}_{\bC,1}^{m}} f_i$ is $\prt{n-c_m} L$-smooth and hence 
    \[
    \begin{aligned}
        & \norm{\sum_{i\notin \tilde{\cI}_{\bC,1}^{m}} \prt{\nabla f_i(x_1^{(t-m)}) - \nabla f_i(x_1^{(t-m-1)})}}^2 \\
        \le &  \prt{n-c_m}^2 L^2 \norm{x_1^{(t-m)} - x_1^{(t-m-1)} }^2.
    \end{aligned}
    \]
    Then, it is implied by choosing $\alpha_m = \prt{1-\frac{\gamma n \mu }{4}}^{d_\bC - m }$ in Lemma \ref{lem:sum_helper} that
    for $T > 2d_\bC$, $\sum_{t=0}^{T}\prt{1-\frac{\gamma n \mu}{4}}^{T-t} \norm{\bPi_\bC\bC^t}_F^2\norm{ \nabla \bF^{(0)}}^2_F \le  2nc_\avg \prt{1-\frac{\gamma n \mu}{4}}^{T-d_\bC} \norm{ \nabla \bF^{(0)} }^2_F,$
    and  noting that $\frac{1}{2}\le 1-\frac{\gamma n \mu}{4} < 1$,
    \begin{equation}
        \label{eq:lem:xk-Ak-2-4}
        \begin{aligned}
             & \sum_{t=0}^{T}\prt{1-\frac{\gamma n \mu}{4}}^{T-t}\sum_{m=0}^{\min\crk{t,d_\bC} - 1} \prt{1-\frac{\gamma n \mu }{4}}^{ m -d_\bC} \cdot \\
            & \qquad \prt{n-c_m}^2 \norm{x_1^{(t-m)} - x_1^{(t-m-1)} }^2\\
            \le & 2 n^2 d_\bC \sum_{t=0}^{T} \prt{1-\frac{\gamma n \mu}{4}}^{T-t-d_\bC} \norm{x_1^{(t+1)} - x_1^{(t)} }^2.
        \end{aligned}
    \end{equation}
    Invoking Lemma \ref{lem:gamma} completes the proof.

\subsection{Proof of Lemma \ref{lem:X_diff}}
\label{pf:lem:X_diff}
Notice that $\Delta \bar{\bX}^{(t)} =   -\gamma\mone \bpi_\bR^\T\bPi_\bC\bY^{(t)}  -\gamma \mone\mone^T \bY^{(t)}.$
    Invoking Equation \eqref{eq:pp22} and noting the fact that $\mone^\T\bY^{(t)} = \mone^\T\bG^{(t)}$, we have
    \begin{equation}
        \label{eq:lem:X_diff-2}
        \begin{aligned}
            \Delta \bar{\bX}^{(t)} = & -\gamma \mone  \sum_{m=1}^{\min\crk{t,d_{\bC}}} \bpi_\bR^\T\bA_{m}\bG^{(t-m)} - \gamma \mone \bpi_\bR^\T \bG^{(t)}. \\
        \end{aligned}
    \end{equation}
    By decomposing $\bG^{(t)} = \bTh^{(t)} + \nabla\bar{\bF}^{(t)} + \nabla \bF(\bar{\bX}^{(t)})$ and in light of Equation \eqref{eq:pp22}, we derive that
    \begin{equation}
        \label{eq:lem:X_diff-2-1}
        \begin{aligned}
            & \Delta \bar{\bX}^{(t)} =  -\gamma \mone  \sum_{m=0}^{\min\crk{t,d_{\bC}}} \bpi_\bR^\T\bA_{m}\brk{\bTh^{(t-m)} + \nabla\bar{\bF}^{(t-m)}} \\
            & -\gamma \mone \bpi_\bR^\T\sum_{m=0}^{\min\crk{t,d_\bC} - 1} \bPi_\bC\bC^{m}\brk{\nabla\bF(\bar{\bX}^{(t-m)}) - \nabla\bF(\bar{\bX}^{(t-m-1)})} \\
            & - \gamma \mone \bpi_\bR^\T \bPi_{\bC}\bC^t \nabla \bF^{(0)} + \gamma \mone \bpi_\bR^\T \bPi_{\bC}\nabla \bF(\bar{\bX}^{(t)}) - \gamma \mone \bpi_\bR^\T \nabla \bF(\bar{\bX}^{(t)})\\
            & :=  \bH_{1}^{(t)} + \bH_{2}^{(t)}  + \bH_{3}^{(t)}  + \bH_{4}^{(t)}  + \bH_{5}^{(t)} ,
        \end{aligned}
    \end{equation}
    where we specify the forms of $\bH_{1}^{(t)} $ to $\bH_{5}^{(t)} $ and provide their upper bounds below. 
    Invoking Lemma \ref{lem:var1} and Corollary \ref{lem:var1}, we have $\expect \norm{\bH_1^{(t)} }_F^2 = \expect \norm{-\gamma \mone  \sum_{m=0}^{\min\crk{t,d_{\bC}}} \bpi_\bR^\T\bA_{m}\bTh^{(t-m)} }_F^2 = \gamma^2 n \expect \norm{ \sum_{m=0}^{\min\crk{t,d_{\bC}}} \bpi_\bR^\T\bA_{m}\bTh^{(t-m)}  }_2^2 \le \gamma^2 n^2 \sigma^2.$\\
    By Assumption \ref{a.smooth}, Lemmas \ref{lem:sum_matrix} and \ref{lem:Ak-barF}, we have 
    \begin{equation}
        \label{eq:lem:X_diff-3-2}
        \begin{aligned}
            & \expect \norm{\bH_2^{(t)} }_F^2 =  \expect \norm{-\gamma \mone  \sum_{m=0}^{\min\crk{t,d_{\bC}}} \bpi_\bR^\T\bA_{m}\nabla\bar{\bF}^{(t-m)}  }_F^2\\
            & =  \gamma^2 n \expect \norm{ \sum_{m=0}^{\min\crk{t,d_{\bC}}} \bpi_\bR^\T\bA_{m}\nabla\bar{\bF}^{(t-m)} }_2^2 ,
        \end{aligned}
    \end{equation}
    and hence $\sum_{t=0}^{T}\expect \norm{\bH_2^{(t)} }_F^2 \le \gamma^2 n^2 L^2 \sum_{t=0}^{T} \expect \norm{\bPi_\bR\bX^{(t)}}_F^2.$
    By Assumption \ref{a.smooth}, Lemmas \ref{lem:Rk_norm} and \ref{lem:sum_matrix}, we have 
    \begin{equation}
        \label{eq:lem:X_diff-3-3}
        \begin{aligned}
            & \expect \norm{\bH_3^{(t)} }_F^2 = \expect \bigg\|-\gamma \mone \bpi_\bR^\T\sum_{m=0}^{\min\crk{t,d_\bC} - 1} \bPi_\bC\bC^{m} \cdot \\
            & \qquad \brk{\nabla\bF(\bar{\bX}^{(t-m)}) - \nabla\bF(\bar{\bX}^{(t-m-1)})} \bigg\|_F^2\\
            \le & 2\gamma^2 n d_\bC L^2 \sum_{m=0}^{\min\crk{t,d_\bC} - 1} \prt{n-c_m}\expect \norm{\bar{\bX}^{t-m} - \bar{\bX}^{t-m-1}}_F^2.
        \end{aligned}
    \end{equation}
    It is obvious that $\expect \norm{\bH_4^{(t)} }_F^2  =  \expect \norm{-\gamma \mone \bpi_\bR^\T \bPi_{\bC}\bC^t \nabla \bF^{(0)} }_F^2 \le  2\gamma^2 n \prt{n-c_t} \norm{\nabla \bF^{(0)} }_F^2$, and $\expect \norm{\bH_5^{(t)} }_F^2 =  \expect \norm{\gamma \mone \brk{\bpi_\bR^\T \bPi_{\bC} - \bpi_\bR^\T}\nabla \bF(\bar{\bX}^{(t)}) }_F^2 = \gamma^2 n^3 \expect \norm{\nabla f(x_1^{(t)})}^2$.
    Thus, taking full expectation on both sides of Equation \eqref{eq:lem:X_diff-2-1}, summing over $t$, implementing the above results on $H_i^{(t)}$ respectively, and invoking Lemma \ref{lem:sum_help}, we obtain
    \begin{equation}
        \label{eq:lem:X_diff-4}
        \begin{aligned}
            &\sum_{t=0}^{T} \expect\norm{\Delta\bar{\bX}^{(t)} }_F^2  \le   5 \gamma^2 n^2 L^2 \sum_{t=0}^{T} \expect \norm{\hat{\bX}^{(t)}}_F^2 \\
            &  + 10 \gamma^2 n^2 d_\bC c_{\avg}L^2 \sum_{t=0}^{T} \expect\norm{\Delta\bar{\bX}^{(t)} }_F^2 + 5\gamma^2 n^2 \sigma^2\prt{T+1} \\
            & + 10 \gamma^2 n^2 c_{\avg}\norm{\nabla \bF^{(0)} }_F^2 + 5\gamma^2 n^3 \sum_{t=0}^{T} \expect\norm{\nabla f(x_1^{(t)})}^2.
        \end{aligned}
    \end{equation}
    For $\gamma \le 1/\prt{6 n\sqrt{d_\bC c_\avg} L } $, we have $10 \gamma^2 n^2 d_\bC c_{\avg}L^2\le \frac{1}{2}$ and obtain the first desired result.

    Now consider the following fact:
    \begin{equation}
        \label{eq:lex:xk-Ak-X-diff-2}
        \begin{aligned}
            & \sum_{t=0}^{T} \prt{1-\frac{\gamma n \mu}{4}}^{T-t} \expect \norm{ \sum_{m=0}^{\min\crk{t,d_\bC}}\bpi_\bR^\T \bA_m  \bTh^{(t-m)} }^2 \\
            \le & \sum_{t=0}^{T} \prt{1-\frac{\gamma n \mu}{4}}^{T-t} n\sigma^2 \le \frac{4n\sigma^2}{\gamma n \mu} = \frac{4\sigma^2}{\gamma \mu}.
        \end{aligned}
    \end{equation}
    From Equation \eqref{eq:lem:X_diff-2-1}, invoking Lemma \ref{lem:Ak-barF} and \ref{lem:xk-Ak-2}, we have
    \begin{equation}
        \label{eq:lex:xk-Ak-X-diff-6}
        \begin{aligned}
            & \sum_{t=0}^{T} \prt{1-\frac{\gamma n \mu}{4}}^{T-t} \expect\prt{\frac{1}{n} - 18 \gamma L } \norm{ \Delta\bar{\bX}^{(t)} }^2_F \\
            & \le \frac{12\gamma \sigma^2}{ \mu} + 18\gamma^2 nc_\avg \prt{1-\frac{\gamma n \mu}{4}}^{T-d_\bC} \norm{ \nabla \bF^{(0)}}^2_F \\
            & + \sum_{t=0}^{T} \prt{1-\frac{\gamma n \mu}{4}}^{T-t} \expect \brk{\frac{3\gamma L}{ d_\bC\kappa } \norm{\hat{\bX}^{(t)}}_F^2 + 9\gamma^2 n^2 \norm{\nabla f(x_1^{(t)})}^2}.
        \end{aligned}
    \end{equation}
    Then, for $\gamma \le 1/\prt{100 n d_\bC \kappa L}$, we have $18 \gamma n L \le 1/2$ and get the second desired result by rearranging the terms.

\subsection{Proof of Lemma \ref{lem:PiX}}
\label{pf:lem:PiX}
With identical initial points, i.e., $\bPi_\bR\bX^{(0)} = \mathbf{0}$, multiplying $\bPi_\bR$ on both sides of Equation \eqref{eq:pp1} and invoking Lemma \ref{lem:Rk_norm} yields
    \begin{equation}
        \label{eq:lem:PiX-1}
        \begin{aligned}
            \hat{\bX}^{(t)}= & -\gamma \sum_{m=0}^{t-1}  \sum_{j=\max\crk{1,t-m-d_\bC}}^{t-m} \bPi_\bR\bR^j\bA_{t-m-j} \bG^{(m)}.\\
        \end{aligned}
    \end{equation}
    Then, by the decomposition of $\bG^{(m)}$, we have
    \begin{equation}
        \label{eq:lem:PiX-2}
        \begin{aligned}
            & \bPi_\bR\bX^{(t)} = -\gamma \sum_{m=\max\crk{0, t-d_\bR - d_\bC + 1}}^{t-1}  \sum_{j=\max\crk{1,t-m-d_\bC}}^{\min\crk{t-m, d_\bR - 1}} \cdot \\
            & \qquad \bPi_\bR\bR^j\bA_{t-m-j} \prt{\bTh^{(m)} + \nabla\bar{\bF}^{(m)} + \nabla \bF(\bar{\bX}^{(m)})}\\
            = &-\gamma \bM_1^{(t)} - \gamma\bM_2^{(t)} - \gamma\bM_3^{(t)}.
        \end{aligned}
    \end{equation}
    where $\bM_1^{(t)}$ to $\bM_3^{(t)}$ are defined as the terms in Equation \eqref{eq:lem:PiX-2} respectively.
    
    Invoking Assumption \ref{a.var}, we have
    \begin{equation}
        \label{eq:lem:PiX-3}
        \begin{aligned}
        & \expect\norm{\bM_1^{(t)}}_F^2 =\expect \norm{\sum_{m}  \sum_{j} \bPi_\bR\bR^j\bA_{t-m-j} \bTh^{(m)}}_F^2 \\
        \le & \sum_{m} \min\crk{d_\bR, d_\bC} \sum_{j}[] 2\prt{n-r_{j}}\cdot2\prt{n-c_{t-m-j-1}}\sigma^2\\
        \le & 4n^2 \min\crk{d_\bR, d_\bC} r_{\avg} c_{\avg} \sigma^2.
        \end{aligned}
    \end{equation}
    Invoking Assumption \ref{a.smooth} and Lemma \ref{lem:var0}, we have
    \begin{equation}
        \label{eq:lem:PiX-4}
        \begin{aligned}
            & \expect\norm{\bM_2^{(t)}}_F^2=  \expect \norm{\sum_{m}  \sum_{j} \bPi_\bR\bR^j\bA_{t-m-j} \nabla\bar{\bF}^{(m)}}_F^2 \\
            \le & \prt{d_\bR + d_\bC - 2} \sum_{m=\max\crk{0, t-d_\bR - d_\bC + 1}}^{t-1} \min\crk{d_\bR, d_\bC}\cdot  \\
            & \sum_{j=\max\crk{1,t-m-d_\bC}}^{\min\crk{t-m, d_\bR - 1}} \norm{\bPi_\bR\bR^j}_2^2\norm{\bA_{t-m-j}}_F^2 L^2\expect\norm{\hat{\bX}^{(m)}}_F^2\\
            \le & 16  d_\bR d_\bC L^2 \sum_{m}\sum_{j} \prt{n-r_j}\prt{n-c_{t-m-j-1}} \expect\norm{\hat{\bX}^{(m)}}_F^2.
        \end{aligned}
    \end{equation}
    It follows from Lemma \ref{lem:sum_help} that
    \begin{equation}
        \label{eq:lem:PiX-4-1}
        \sum_{t=0}^{T}\expect\norm{\bM_2^{(t)}}_F^2 \le 16 n^2  d_\bR d_\bC r_{\avg} c_{\avg} L^2 \sum_{t=0}^{T} \expect\norm{\hat{\bX}^{(t)}}_F^2.
    \end{equation}
    For $\bM_3^{(t)}$, in light of Equation \eqref{eq:lem:PiX-1}, we have
    \begin{equation}
        \label{eq:lem:PiX-5}
        \begin{aligned}
            & \bM_3^{(t)} :=  \bM_{3,1}^{(t)} + \bM_{3,2}^{(t)} + \bM_{3,3}^{(t)} \\
            & =   \sum_{m=0}^{t-2} \sum_{j=1}^{t-m-1}\bPi_\bR\bR^{j}\bPi_\bC\bC^{t-m-1-j} \brk{\nabla \bF(\bar{\bX}^{(m+1)} )- \nabla \bF(\bar{\bX}^{(m)})}  \\
            & + \sum_{j=1}^{t}\bPi_\bR\bR^{j}\bPi_\bC\bC^{t-j} \nabla \bF^{(0)}   + \sum_{j=1}^{t}\bPi_\bR\bR^{j} \bpi_\bC \mone^\T \nabla \bF(\bar{\bX}^{(t-j)}),\\
        \end{aligned}
    \end{equation}
    where $\bM_{3,1}^{(t)}$, $\bM_{3,2}^{(t)}$ and $\bM_{3,3}^{(t)}$ are defined as the terms in Equation \eqref{eq:lem:PiX-5} respectively.
    
    Then, by Lemma \ref{lem:Rk_norm} and \ref{lem:sum_matrix}, we have
    \begin{equation}
        \label{eq:lem:PiX-5-1}
        \begin{aligned}
            & \expect\norm{\bM_{3,1}^{(t)}}_F^2 =  \expect \norm{\sum_{j=1}^{t}\bPi_\bR\bR^{j}\bPi_\bC\bC^{t-j} \nabla \bF^{(0)} }_F^2\\
            \le & 4\min\crk{d_\bR,d_\bC} \sum_{j=\max\crk{1,t-d_\bC + 1}}^{\min\crk{t,d_\bR - 1}} \prt{n-r_j}\prt{n-c_{t-j}}\norm{\nabla \bF^{(0)}}_F^2.
        \end{aligned}
    \end{equation}
    After summing over $t$, we obtain
    \begin{equation}
        \label{eq:lem:PiX-5-2}
        \begin{aligned}
        & \sum_{t=0}^{T}\expect\norm{\bM_{3,1}^{(t)}}_F^2 
        \le  4 n^2\min\crk{d_\bR, d_\bC} r_\avg c_\avg \norm{\nabla \bF^{(0)}}_F^2.
        \end{aligned}
    \end{equation}
    By Assumption \ref{a.smooth}, Lemmas \ref{lem:Rk_norm} and \ref{lem:sum_matrix}, we have
    \[
        \begin{aligned}
            & \expect\norm{\bM_{3,2}^{(t)}}_F^2  
            \le  8 d_\bR d_\bC L^2 \sum_{m=\max\crk{0,t-d_\bR-d_\bC +1 }}^{t-2} \cdot \\
            & \sum_{j=\max\crk{1,t-m-d_\bC}}^{\min\crk{t-m-1, d_\bR - 1}}\prt{n-r_j}\prt{n-c_{t-m-1-j}}\expect\norm{\Delta \bar{\bX}^{(m)} }_F^2.
        \end{aligned}
    \]
    Then, summing over $t$ and invoking Lemma \ref{lem:sum_help}, we get 
    \begin{equation}
        \label{eq:lem:PiX-5-3}
        \begin{aligned}
            \sum_{t=0}^{T}\expect\norm{\bM_{3,2}^{(t)}}_F^2 \le 8 n^2 d_\bR d_\bC r_{\avg} c_{\avg} L^2 \sum_{t=0}^{T} \expect\norm{\Delta\bar{\bX}^{(t)}}_F^2.
        \end{aligned}
    \end{equation}
    For $\bM_{3,3}^{(t)}$, we have
    \begin{equation}
        \label{eq:lem:PiX-5-4}
        \begin{aligned}
            & \expect\norm{\bM_{3,3}^{(t)}}_F^2 
            =  \expect\norm{\sum_{j=1}^{\min\crk{t,d_\bR - 1}}\bPi_\bR\bR^{j} \bpi_\bC \mone^\T \nabla \bF(\bar{\bX}^{(t-j)})}_F^2\\
            \le & d_\bR \sum_{j=1}^{\min\crk{t,d_\bR - 1}}\norm{\bPi_\bR \bR^j}_2^2\norm{\bpi_\bC}_2^2n^2\norm{\frac{1}{n}\mone^\T\nabla \bF(\bar{\bX}^{(t-j)}) }_F^2\\
            \le & 2 n^2 d_\bR \sum_{j=1}^{\min\crk{t,d_\bR - 1}} \prt{n-r_j} \norm{\nabla f(x_1^{(t-j)})}^2.
        \end{aligned}
    \end{equation}
    Summing over $t$ and invoking Lemma \ref{lem:sum_help}, we have
    \begin{equation}
        \label{eq:lem:PiX-5-5}
        \sum_{t=0}^{T}\expect\norm{\bM_{3,3}^{(t)}}_F^2 \le 2 n^3 d_\bR r_{\avg} \sum_{t=0}^{T} \expect\norm{\nabla f(x_1^{(t)})}^2.
    \end{equation}
    Back to Equation \eqref{eq:lem:PiX-2}, after taking full expectation on both sides, summing over $t$ and invoking Equation \eqref{eq:lem:PiX-3} to Equation \eqref{eq:lem:PiX-5-5}, we obtain $\sum_{t=0}^{T}\expect\norm{\hat{\bX}^{(t)}}_F^2  \le 3\gamma^2 \sum_{t=0}^{T} \brk{\expect\norm{\bM_1^{(t)}}_F^2 + \expect\norm{\bM_2^{(t)}}_F^2} + 9\gamma^2 \sum_{t=0}^{T}\brk{\expect \norm{\bM_{3,1}^{(t)}}_F^2+ \expect \norm{\bM_{3,2}^{(t)}}_F^2+ \expect \norm{\bM_{3,3}^{(t)}}_F^2}$, or
    \begin{equation}
        \label{eq:lem:PiX-6}
        \begin{aligned}
            & \sum_{t=0}^{T}\expect\norm{\hat{\bX}^{(t)}}_F^2  \le 12 \gamma^2 n^2 \min\crk{d_\bR, d_\bC} r_\avg c_\avg \sigma^2\prt{T+1}\\
            & + 48\gamma^2 n^2 d_\bR d_\bC r_\avg c_\avg L^2 \sum_{t=0}^{T}\expect\norm{\hat{\bX}^{(t)}}_F^2\\
            & + 72\gamma^2 n^2 d_\bR d_\bC r_\avg c_\avg L^2 \sum_{t=0}^{T} \expect\norm{\Delta\bar{\bX}^{(t)} }_F^2\\
            & + 36 \gamma^2 n^2 \min\crk{d_\bR, d_\bC} r_\avg c_\avg \norm{\nabla \bF^{(0)}}_F^2 \\
            & + 18 \gamma^2 n^3 d_\bR r_\avg \sum_{t=0}^{T} \expect\norm{\nabla f(x_1^{(t)})}^2.
        \end{aligned}
    \end{equation}
    Implementing Lemma \ref{lem:X_diff} into Equation \eqref{eq:lem:PiX-6}, we derive that, 
    for $\gamma \le 1/\prt{100n\sqrt{d_\bR d_\bC r_\avg c_\avg} L}$,
    \[
    \begin{aligned}
        & \sum_{t=0}^{T}\expect\norm{\hat{\bX}^{(t)}}_F^2  \le 20 \gamma^2 n^2 \min\crk{d_\bR, d_\bC} r_\avg c_\avg \sigma^2\prt{T+1} \\
        & \quad + 20 \gamma^2 n^3 d_\bR r_\avg \sum_{t=0}^{T} \expect\norm{\nabla f(x_1^{(t)})}^2 + \frac{1}{2} \sum_{t=0}^{T}\expect\norm{\hat{\bX}^{(t)}}_F^2\\
        & \quad  + 40 \gamma^2 n^2 \min\crk{d_\bR, d_\bC} r_\avg c_\avg \norm{\nabla \bF^{(0)}}_F^2,
    \end{aligned}
    \]
    which implies the first desired inequality. 
    
    Furthermore, with $\bM_1^{(t)}$ to $\bM_3^{(t)}$ defined as the terms in Equation \eqref{eq:lem:X_diff-2-1}, invoking Assumption \ref{a.var} and Equation \eqref{eq:lem:PiX-3}, we have
    \begin{equation}
        \label{eq:lem:xk-Ak-Pi-X-2}
        \begin{aligned}
            & \sum_{t=0}^{T}\prt{1-\frac{\gamma n \mu}{4}}^{T-t}\expect \norm{\bM_1^{(t)}}_F^2 \le  \sum_{t=0}^{T}\prt{1-\frac{\gamma n \mu}{4}}^{T-t} 4 n^2 \cdot\\
            & \min\crk{d_\bC,d_\bC} r_\avg c_\avg \sigma^2 \le \frac{16 n \min\crk{d_\bC,d_\bC} r_\avg c_\avg \sigma^2}{\gamma \mu}.
        \end{aligned}
    \end{equation}
    By choosing $\alpha_m = \prt{1-\prt{\gamma n \mu}/4}^{t - m}$ in Lemma \ref{lem:sum_helper}, we have
    \begin{equation}
        \label{eq:lem:xk-Ak-Pi-X-3}
        \begin{aligned}
            &\norm{M_2^{(t)}}_F^2 
            \le  \sum_{m=\max\crk{0,t-d_\bR - d_\bC + 1}}^{t-1}\prt{1-\frac{\gamma n \mu}{4}}^{ m - t }  \cdot \\
            & \quad \norm{ \sum_{j = \max\crk{1, t-m-d_\bC}}^{\min\crk{t-m, d_\bC - 1}} \bPi_\bR \bR^{j} \bA_{t-m-j} \nabla \bar{\bF}^{(m)}   }_F^2 \\
            & \le  \sum_{m=\max\crk{0,t-d_\bR - d_\bC + 1}}^{t-1} \prt{1-\frac{\gamma n \mu}{4}}^{ m - t } \sum_{j = \max\crk{1, t-m-d_\bC}}^{\min\crk{t-m, d_\bC - 1}}\cdot \\
            &  \quad \prt{1-\frac{\gamma n \mu}{4}}^{ -j } \norm{ \bPi_\bR \bR^{j} }_2^2 \norm{ \bA_{t-m-j}  }_F^2 \norm{ \nabla \bar{\bF}^{(m)}}_F^2\\
            & \le  4n^2 \min\crk{d_\bC, d_\bC} \prt{1-\frac{\gamma n \mu}{4}}^{ -d_\bC } \sum_{m=\max\crk{0,t-d_\bR - d_\bC + 1}}^{t-1} \cdot\\
            &\quad \prt{1-\frac{\gamma n \mu}{4}}^{ m - t }    \norm{ \nabla \bar{\bF}^{(m)}}_F^2.\\
        \end{aligned}
    \end{equation}
    Then, it holds that 
    \begin{equation}
        \label{eq:lem:xk-Ak-Pi-X-4}
        \begin{aligned}
            & \sum_{t=0}^{T}\prt{1-\frac{\gamma n \mu}{4}}^{T-t}\expect \norm{\bM_2^{(t)}}_F^2 \le  4n^2 \min\crk{d_\bC, d_\bC}  \cdot\\
            & \prt{1-\frac{\gamma n \mu}{4}}^{ -d_\bC } \sum_{t=0}^{T}  \sum_{m=a_t}^{t-1} \prt{1-\frac{\gamma n \mu}{4}}^{ T - 2t + m  }    \expect\norm{ \nabla \bar{\bF}^{(m)}}_F^2,
        \end{aligned}
    \end{equation}
    where $a_t = \max\crk{0,t-d_\bR - d_\bC + 1}$. Consequently,  
    \begin{equation}
        \label{eq:lem:xk-Ak-Pi-X-5}
        \begin{aligned}
            & \sum_{m=\max\crk{0,t-d_\bR - d_\bC + 1}}^{t-1} \prt{1-\frac{\gamma n \mu}{4}}^{ T - 2t + m }    \norm{ \nabla \bar{\bF}^{(m)}}_F^2  \\
            \le &  \prt{ 1-\frac{\gamma n \mu}{4}}^{T-t}  \prt{ 1-\frac{\gamma n \mu}{4}}^{-(d_\bR + d_\bC)} \sum_{m=a_t}^{t-1}  \norm{ \nabla \bar{\bF}^{(m)}}_F^2.
        \end{aligned}
    \end{equation}
    Thus, implementing with Equation \eqref{eq:lem:xk-Ak-Pi-X-4}, we have
    \begin{equation}
        \label{eq:lem:xk-Ak-Pi-X-6}
        \begin{aligned}
            & \sum_{t=0}^{T}\prt{1-\frac{\gamma n \mu}{4}}^{T-t}\expect \norm{\bM_2^{(t)}} \le 16 n^2 \min\crk{d_\bC, d_\bC} \prt{d_\bR + d_\bC}  \cdot\\
            & L^2\prt{1-\frac{\gamma n \mu}{4}}^{-3(d_\bR + d_\bC) } \sum_{t=0}^{T} \prt{1-\frac{\gamma n \mu}{4}}^{T-t} \norm{\hat{\bX}^{(t)}}_F^2\\
            &\le \frac{3n \min\crk{d_\bC, d_\bC} L }{\gamma \kappa } \sum_{t=0}^{T} \prt{1-\frac{\gamma n \mu}{4}}^{T-t} \norm{\hat{\bX}^{(t)}}_F^2,
        \end{aligned}
    \end{equation}
    where we use Lemma \ref{lem:gamma} with $d = 3(d_\bC+d_\bC)\kappa $ and $\gamma \le \frac{1}{100 n\max\crk{d_\bR,d_\bC}\kappa L}$ in the last inequality.
    
    For $\bM_3^{(t)}$, in light of Equation \eqref{eq:lem:PiX-5}, we decompose $\bM_3^{(t)}$ into $\bM_{3,1}^{(t)}$, $ \bM_{3,2}^{(t)} $, and $ \bM_{3,3}^{(t)}$.
    For $\bM_{3,1}^{(t)}$, by Equation \eqref{eq:lem:PiX-5-1}, we have, with $T > 2(d_\bR + d_\bC)$, that
    \begin{equation}
        \label{eq:lem:xk-Ak-Pi-X-9}
        \begin{aligned}
            & \sum_{t=0}^{T}\prt{1-\frac{\gamma n \mu}{4}}^{T-t}\expect\norm{\bM_{3,1}^{(t)}}_F^2 
            \le  4n^2 \prt{\min\crk{d_\bR,d_\bC}}^2  \cdot \\
            & \qquad \prt{d_\bR + d_\bC}\prt{1-\frac{\gamma n \mu}{4}}^{T-d_\bR -d_\bC} \norm{ \nabla \bF^{(0)}}_F^2.
        \end{aligned} 
    \end{equation}
    For $\bM_{3,2}^{(t)}$, invoking Equation \eqref{eq:lem:PiX-5-2} and Lemma \ref{lem:sum_helper} with $\alpha_m = \prt{1-\frac{\gamma n \mu}{4}}^{t-m}$, $\alpha_j = \prt{1-\frac{\gamma n \mu}{4}}^{d_\bR - j}$ with respect to different summations, we get
    \begin{equation}
        \label{eq:lem:xk-Ak-Pi-X-10}
        \begin{aligned}
            &\norm{\bM_{3,2}^{(t)}}_F^2 \le  \sum_{m=\max\crk{0,t-d_\bR - d_\bC + 1}}^{t-2} \prt{1-\frac{\gamma n \mu}{4}}^{m - t}\cdot \\
            & \quad \sum_{j=\max\crk{1, t-m-d_\bC}}^{\min\crk{t-m-1, d_\bR -1}} \prt{1-\frac{\gamma n \mu}{4}}^{ j - d_\bR}\norm{\bPi_\bR\bR^{j}}_F^2\cdot \\
            & \quad \norm{\bPi_\bC\bC^{t-m-1-j}}_F^2\norm{\nabla \bF(\bar{\bX}^{(m+1)} )- \nabla \bF(\bar{\bX}^{(m)}) }_F^2 \\
            \le &   4n^3 d_\bR L^2 \prt{1-\frac{\gamma n \mu}{4}}^{-(2d_\bR + d_\bC)} \sum_{m=a_t}^{t-2} \norm{ x_1^{(m+1)} - x_1^{(m)}}^2,
        \end{aligned}
    \end{equation}
    where $a_t = \max\crk{0,t-d_\bR - d_\bC + 1}$. Then, it implies that for $\gamma \le \frac{1}{100 n \max\crk{d_\bR, d_\bC} L}$,
    \begin{equation}
        \label{eq:lem:xk-Ak-Pi-X-10-1}
        \begin{aligned}
            &\sum_{t=0}^{T}\prt{1-\frac{\gamma n \mu}{4}}^{T-t} \expect\norm{\bM_{3,2}^{(t)}}_F^2 \\
            \le & \frac{4n^2 d_\bR  L}{\gamma \kappa } \sum_{t=0}^{T}\prt{1-\frac{\gamma n \mu}{4}}^{T-t} \expect\norm{ x_1^{(t+1)} - x_1^{(t)}}^2.
        \end{aligned}
    \end{equation}
    For $\bM_{3,3}^{(t)}$, invoking Lemma \ref{lem:sum_helper} with $\alpha_j = \prt{1 - \frac{\gamma n \mu}{4}}^{d_\bC - j}$, we have 
    \begin{equation}
        \label{eq:lem:xk-Ak-Pi-X-11}
        \begin{aligned}
            & \sum_{t=0}^{T}\prt{1-\frac{\gamma n \mu}{4}}^{T-t} \norm{\bM_{3,3}^{(t)}}_F^2 \le \sum_{t=0}^{T}\prt{1-\frac{\gamma n \mu}{4}}^{T-t} \cdot\\
            & \sum_{j=1}^{\min\crk{t,d_\bR - 1}} \prt{1-\frac{\gamma n \mu}{4}}^{j-d_\bC} \norm{ \bPi_\bR \bR^{j} \bpi_\bC }^2 n^2 \norm{\nabla f(x_1^{(t-j)})}^2 \\
            \le & n^3 d_\bR  \sum_{t=0}^{T} \prt{1-\frac{\gamma n \mu}{4}}^{T-t - d_\bC}\norm{\nabla f(x_1^{(t)})}^2.
        \end{aligned}
    \end{equation}
    Finally, by taking squared F-norm on both sides of Equation \eqref{eq:lem:PiX-2}, multiplying $\prt{1-\frac{\gamma n \mu}{4}}^{T-t}$ and then summing over $t$, we have, after implementing the results on $\sum_{t=0}^{T}\prt{1-\frac{\gamma n \mu}{4}}^{T-t} \expect \norm{\bM^{(t)}}_F^2$, 
    for $\gamma \le \frac{1}{100 n\max\crk{d_\bR,d_\bC}\kappa L}$,
    \begin{equation}
        \label{eq:lem:xk-Ak-Pi-X-13}
        \begin{aligned}
            & \sum_{t=0}^{T}\prt{1-\frac{\gamma n \mu}{4}}^{T-t} \expect \norm{\hat{\bX}^{(t)}}_F^2 \le  \frac{48 \gamma n \min\crk{d_\bC,d_\bC} r_\avg c_\avg \sigma^2}{ \mu} \\
            & + 9\gamma n \min\crk{d_\bC, d_\bC} \frac{L}{\kappa} \sum_{t=0}^{T} \prt{1-\frac{\gamma n \mu}{4}}^{T-t} \norm{\hat{\bX}^{(t)}}_F^2 + 36\gamma^2 n^2\cdot\\
            & \  \prt{\min\crk{d_\bR,d_\bC}}^2\prt{d_\bR + d_\bC} \prt{1-\frac{\gamma n \mu}{4}}^{T-d_\bR -d_\bC} \norm{ \nabla \bF^{(0)} }_F^2\\
            & + 36\gamma n^2 d_\bR \frac{L}{\kappa} \sum_{t=0}^{T}\prt{1-\frac{\gamma n \mu}{4}}^{T-t} \norm{ x_1^{(t+1)} - x_1^{(t)}}^2 \\
            & + 9\gamma^2 n^3 d_\bR  \sum_{t=0}^{T} \prt{1-\frac{\gamma n \mu}{4}}^{T-t - d_\bC}\norm{\nabla f(x_1^{(t)})}^2.
        \end{aligned}
    \end{equation}
    Plugging Lemma \ref{lem:X_diff} into Equation \eqref{eq:lem:xk-Ak-Pi-X-13}, we get the desired result for $\gamma \le \frac{1}{100 n\max\crk{d_\bR,d_\bC}\kappa L}$,
    noting the facts
    \begin{equation}
        \label{eq:lem:xk-Ak-Pi-X-15}
        \begin{aligned}
            & 9\gamma^2 n^3 d_\bR  + 800\gamma^3 n^4 d_\bR L \le 17\gamma^2 n^3 d_\bR,  \\
            & 100\gamma n c_\avg   \le  \prt{\min\crk{d_\bR,d_\bC}}^2\prt{d_\bR + d_\bC} \prt{1-\frac{\gamma n \mu}{4}}^{-d_\bR }, \\
            & 9\gamma n \min\crk{d_\bC, d_\bC} L + 216\gamma^2 n^2 \frac{d_\bR}{d_\bC}L^2 \le \frac{1}{2}.
        \end{aligned}
    \end{equation}

\subsection{Proof of Lemma \ref{lem:independ1}}
\label{pf:lem:independ1}
Denote $g_i^{(t)}:=g_i(x_i^{(t)}, \xi_i^{(t)})$ for simplicity.
    Notice that
    \[
    \begin{aligned}
        & x_1^{(t)}  = x_{1}^{(t-1)} - \gamma y_{1}^{(t-1)} = x_{1}^{(t-1)} \\
        & \qquad - \gamma \sum_{i\in \cI_{\bC,1}^{1}} y_i^{(t-2)} - \gamma g_1^{(t-1)} + \gamma g_1^{(t-2)}.
    \end{aligned}
    \]
     Therefore, $x_1^{(t)}$ not a function of $\xi_i^{(t-1)}$ for $i\notin \cI_{\bC,1}^{0}$, implying that $x_1^{(t)}$ is independent with $\xi_i^{(t-1)}$ where $i\notin \cI_{\bC,1}^{0}$.  Since $\xi_i^{(t-1)}$ are randomly sampled, we have $\xi_i^{(t-1)}$, $i\notin \cI_{\bC,1}^{0}$, is independent with $(x_1^{(t)}, x_j^{(t-1)})$, $j\in \cI_{\bC,1}^{i}$.
    
    Then, we iterate the above equation and get
    \begin{equation}
        \label{eq:lem:independ1-1}
        \begin{aligned}
            & x_1^{(t)}  = x_{1}^{(t-1)} - \gamma g_1^{(t-1)} + \gamma g_1^{(t-2)} \\
            & - \gamma \sum_{i\in \cI_{\bC,1}^{1}} \brk{ \sum_{j\in \cI_{\bC,i}^{1}} y_j^{(t-3)} + g_i^{(t-2)} - g_i^{(t-3)}  } =  x_{1}^{(t-1)} - \gamma g_1^{(t-1)}\\
            &  - \gamma  \sum_{i\in \cI_{\bC,1}^1,i\ne 1}  g_i^{(t-2)}  -\gamma \sum_{i\in \tilde{\cI}_{\bC,1}^2} y_i^{(t-3)} + \gamma \sum_{i\in \tilde{\cI}_{\bC,1}^1} g_i^{(t-3)} .
        \end{aligned}
    \end{equation}
     Noting that $x_1^{(t-1)}$ is not a function of $\xi_i^{(t-2)}$ and $x_i^{(t-2)}$, we have $i\notin \tilde{\cI}_{\bC,1}^{0}$, $x_1^{(t)}$ is not a function of $\xi_i^{(t-2)}$, $i\notin \tilde{\cI}_{\bC,1}^{1}$. 
    
     By iterating the above procedure, we get the conclusion  that $x^{(t)}_1$ is not a function of $\xi_i^{(t-j)}$ and $x_i^{(t-j)}$, $i\notin \tilde{\cI}_{\bC,1}^{j-1}$, especially for those $i\in \cI_{\bC,1}^{j}$, $j\in \crk{1,2,\cdots, d_\bC}$. 
    
    Consequently, by choosing $X = x_1^{(t)}$, $Y = \crk{x_j^{(t-i)}, j\in \cI_{\bC,1}^{i}}$ and $Z = \crk{ \xi_i^{(t-i)}, j\in \cI_{\bC,1}^{i}}$ in Lemma \ref{lem:independ_help}, we have $Z$ and $(X,Y)$ are independent, and thus
    \begin{equation}
        \label{eq:lem:independ1-2}
        \begin{aligned}
            & \expect\ip{ \nabla f(x_1^{(t)})}{ \bfe_{\cI_{\cC,1}^{i}}^T\bTh^{(t-i)} } = 0 ,
        \end{aligned}
    \end{equation}
    where $h(X) = \nabla f(X)$ and $g(Y,Z) = \bfe_{\cI_{\cC,1}^{i}}^T\bTh^{(t-i)}$ satisfy $\expect{[g(Y,Z)|Y]} = 0$ by Assumption \ref{a.var}.
    
    Furthermore, by choosing $h(X):= X - x^*$ rather than $\nabla f(X)$, we get
    \[
    \expect\ip{ x_1^{(t)} - x^*}{ \bfe_{\cI_{\bC,1}^{i}}^\T\bTh^{(t-i)} } = 0.
    \]

\subsection{Proof of Lemma \ref{lem:smooth}}

\label{pf:lem:smooth}
By Assumption \ref{a.smooth}, the function $f:=\frac{1}{n}f_i$ is $L$-smooth. Then, the descent lemma holds:
    \begin{equation}
        \label{eq:lem:smooth-1}
        \begin{aligned}
            \expect f(x_1^{(t+1)}) \le & \expect f(x_1^{(t)}) + \expect \ip{\nabla f(x_1^{(t)})}{ x_1^{(t+1)} - x_1^{(t)}} \\
            & + \frac{L}{2}\expect \norm{x_1^{(t+1)} - x_1^{(t)}}^2.
        \end{aligned}
    \end{equation}
    For the inner product term, we start with the decomposition of $x_1^{(t+1)} - x_1^{(t)}$, i.e.,
    \begin{equation}
        \label{eq:lem:smooth-2}
        x_1^{(t+1)} - x_1^{(t)} = - \gamma \bpi_\bR^\T \bY^{(t)} = - \gamma \bpi_\bR^\T \bPi_\bC \bY^{(t)} - \gamma \mone^\T \bY^{(t)}.
    \end{equation}
    In light of Equation \eqref{eq:pp22}, $\bpi_\bR^\T \bPi_\bC \bY^{(t)} = \bpi_\bR^\T \sum_{m=0}^{\min\crk{t,d_\bC}}\bA_m \bG^{(t-m)} .$
    
    It results from Lemma \ref{lem:independ1} and the decomposition of $\bG^{(t)} = \bTh^{(t)} + \nabla \bar{\bF}^{(t)} + \nabla \bF(\bX^{(t)})$ that 
    \begin{equation}
        \label{eq:lem:smooth-2-2}
        \begin{aligned}
            & \expect \ip{\nabla f(x_1^{(t)})}{\bpi_\bR^\T \bPi_\bC \bY^{(t)}} := \bP_1^{(t)} + \bP_2^{(t)} \\
            = & \expect \ip{\nabla f(x_1^{(t)})}{ \bpi_\bR^\T\sum_{m=0}^{\min\crk{t,d_\bC}}\bA_m \nabla \bar{\bF}^{(t-m)}}\\
            & + \expect \langle \nabla f(x_1^{(t)}), \bPi_\bC\bC^{t} \nabla \bF^{(0)} + \sum_{m=0}^{\min\crk{t,d_\bC} - 1} \cdot  \\
            & \quad  \bpi_\bR^\T\bPi_\bC\bC^m \brk{\nabla \bF(\bar{\bX}^{(t-m)}) - \nabla \bF(\bar{\bX}^{(t-m-1)})} \rangle\\
        \end{aligned}
    \end{equation}
    where $\bP_1^{(t)}$ and $\bP_2^{(t)}$ are defined as the terms in Equation \eqref{eq:lem:smooth-2-2} respectively.
    
    Then, we have, $-\expect  \bP_1^{(t)} \le  \frac{n}{4} \expect \norm{\nabla f(x_1^{(t)})}^2 + \frac{2}{n}\expect\norm{\bpi_\bR^\T\sum_{m=0}^{\min\crk{t,d_\bC}}\bA_m \nabla \bar{\bF}^{(t-m)}}^2.$
    
    Invoking Lemmas \ref{lem:sum_help} and \ref{lem:Ak-barF}, we have
    \begin{equation}
        \label{eq:lem:smooth-3-1}
        -\gamma \sum_{t=0}^{T}\expect  \bP_1^{(t)} \le \frac{n}{4}\gamma  \sum_{t=0}^{T} \expect \norm{\nabla f(x_1^{(t)})}^2 + 2 \gamma  L^2 \sum_{t=0}^{T} \expect \norm{\hat{\bX}^{(t)}}^2_F.
    \end{equation}
    Besides, 
    \begin{equation}
        \label{eq:lem:smooth-4}
        \begin{aligned}
            & - \gamma \expect  \bP_2^{(t)} 
            \le \frac{ n\gamma}{4} \expect \norm{\nabla f(x_1^{(t)})}^2 + \frac{2\gamma}{n}\norm{\bPi_\bC\bC^{t}}_2^2 \norm{\nabla \bF^{(0)}}^2_F \\
            &  + \frac{2\gamma}{n} d_\bC \sum_{m=0}^{\min\crk{t,d_\bC} - 1} \expect \norm{\bPi_\bC\bC^m}_2^2 L^2 \norm{\Delta\bar{\bX}^{(t-m-1)} }_F^2.\\
        \end{aligned}
    \end{equation}
    Invoking Lemmas \ref{lem:Rk_norm} and \ref{lem:sum_help}, we have
    \begin{equation}
        \label{eq:lem:smooth-4-1}
        \begin{aligned}
        & -\gamma \sum_{t=0}^{T}\expect  \bP_2^{(t)} \le  \frac{\gamma n}{4} \sum_{t=0}^{T} \expect \norm{\nabla f(x_1^{(t)})}^2 + 4\gamma c_\avg \norm{\nabla \bF^{(0)}}^2_F \\
        &\quad  + 4 \gamma d_\bC c_\avg L^2 \sum_{t=0}^{T} \expect \norm{\Delta\bar{\bX}^{(t)} }_F^2 .
        \end{aligned}
    \end{equation}
    Furthermore, from Equation \eqref{eq:lem:smooth-2}, we have, by Lemma \ref{lem:var},
    \begin{equation}
        \label{eq:lem:smooth-5}
        \begin{aligned}
            & \expect \ip{\nabla f(x_1^{(t)}) }{ -\gamma \mone^\T \bY^{(t)} } \\
            = &  -\gamma \expect \ip{\nabla f(x_1^{(t)}) }{ \mone^\T \brk{\bTh^{(t)} + \nabla \bar{\bF}^{(t)} + \nabla \bF(\bar{\bX}^{(t)})} } \\
            = & -\gamma \expect \ip{\nabla f(x_1^{(t)}) }{\mone^\T\nabla \bar{\bF}^{(t)}} - \gamma n \expect \norm{\nabla f(x_1^{(t)})}^2 \\
            \le &  \frac{\gamma n}{4} \expect \norm{\nabla f(x_1^{(t)})}^2 + \gamma L^2 \expect \norm{\bPi_\bR\bX^{(t)}}_F^2 - \gamma n \expect \norm{\nabla f(x_1^{(t)})}^2 .
        \end{aligned}
    \end{equation}
    Back to Equation \eqref{eq:lem:smooth-1}, after summing over $t$ and invoking Equation \eqref{eq:lem:smooth-3-1} to Equation \eqref{eq:lem:smooth-5}, we obtain
    \begin{equation}
        \label{eq:lem:smooth-6}
        \begin{aligned}
            & -\Delta_f \le  -\gamma \sum_{t=0}^{T}\expect  \bP_1^{(t)} -\gamma \sum_{t=0}^{T}\expect  \bP_2^{(t)} \\
            & + \expect \ip{\nabla f(x_1^{(t)}) }{ -\gamma \mone^\T \bY^{(t)} } + \frac{L}{2n}\sum_{t=0}^{T}\expect \norm{\Delta\bar{\bX}^{(t)} }_F^2\\
            & \le   -\frac{\gamma n}{4} \sum_{t=0}^{T}\expect \norm{\nabla f(x_1^{(t)})}^2 + 5 \gamma L^2 \sum_{t=0}^{T} \expect \norm{\hat{\bX}^{(t)}}_F^2\\
            & + \brk{ 4 \gamma d_\bC c_\avg L^2 + \frac{L}{2n} }\sum_{t=0}^{T} \expect \norm{\Delta\bar{\bX}^{(t)}}_F^2 + 4\gamma c_\avg \norm{\nabla \bF^{(0)}}^2_F,
        \end{aligned}
    \end{equation}
    where $\Delta_f:=f(x^{(0)}) - f^*$.
    
    Implementing Lemma \ref{lem:X_diff} into Equation \eqref{eq:lem:smooth-6}, we derive that, for $\gamma \le 1/\prt{100n\sqrt{d_\bR d_\bC r_\avg c_\avg} L} \prt{\le 1/\prt{6n\sqrt{d_\bC c_\avg}L}}$, 
    \begin{equation}
        \label{eq:lem:smooth-7}
        \begin{aligned}
            & -\Delta_f \le  5\gamma^2 n L \sigma^2\prt{T+1} + 40\gamma^3 n^2 d_\bC c_\avg L^2 \sigma^2\prt{T+1} \\
            &  + \brk{5 \gamma  L^2 + 80\gamma^3 n^2 d_\bC c_\avg L^4 + 10 \gamma^2 n  L^3}\sum_{t=0}^{T} \expect \norm{\bPi_\bR\bX^{(t)}}_F^2\\
            &  + \brk{-\frac{\gamma n}{4} + 5\gamma^2 n^2 L + 40 \gamma^3 n^3 d_\bC c_\avg L^2} \sum_{t=0}^{T}\expect \norm{\nabla f(x_1^{(t)})}^2 \\
            &  + \brk{4\gamma n c_\avg + 80\gamma^3 n^3 d_\bC c_\avg^2 L^2 + 10 \gamma^2 n^2 c_\avg L}\frac{1}{n}\norm{\nabla \bF^{(0)}}^2_F .
        \end{aligned}
    \end{equation}
    Implementing Lemma \ref{lem:PiX} into Equation \eqref{eq:lem:smooth-7}, we derive that, for $\gamma \le 1/\prt{100n\sqrt{d_\bR d_\bC r_\avg c_\avg} L}$,
    \[
    5 \gamma L^2 + 80\gamma^3 n^2 d_\bC c_\avg L^4 + 10 \gamma^2 n  L^3 \le 6 \gamma  L^2,
    \]
    and 
    hence
    \begin{equation}
        \label{eq:lem:smooth-9}
        \begin{aligned}
            & -\Delta_f \le 5\gamma^2 n L \sigma^2\prt{T+1} + 300d_\bC r_\avg c_\avg \gamma^3 n^2 L^2 \sigma^2 \prt{T+1} \\
            & \quad -\frac{\gamma n}{8} \sum_{t=0}^{T}\expect \norm{\nabla f(x_1^{(t)})}^2 + 6 \gamma n c_\avg \frac{1}{n}\norm{\nabla \bF^{(0)}}^2_F.
        \end{aligned}
    \end{equation}
    By re-arranging the terms and dividing both sides by $T+1$ on Equation \eqref{eq:lem:smooth-9}, we obtain
    \[
    \begin{aligned}
        &\frac{1}{T+1} \sum_{t=0}^{T}\expect \norm{\nabla f(x_1^{(t)})}^2 \le \frac{8\Delta_f}{\gamma n \prt{T+1}} + 40 \gamma L \sigma^2 \\
        & \quad + 2400 nd_\bC r_\avg c_\avg \gamma^2  L^2 \sigma^2 +  \frac{48 c_\avg}{n\prt{T+1}}\norm{\nabla \bF^{(0)} }^2_F.
    \end{aligned}
    \]
    which completes the proof.

\subsection{Proof of Lemma \ref{lem:convex}}
\label{pf:lem:convex}
We start with analyzing the behavior of \smallmath{\norm{x_1^{(t)} - x^*}^2}. It holds that
    \begin{equation}
        \label{eq:lem:convex-1}
        \begin{aligned}
        & \left\|x_1^{(t+1)} - x^* \right\|^2 = \left\|x_1^{(t)} - x^* \right\|^2 \\
        & \quad + 2 \ip{ x_1^{(t)} - x^*}{ x_1^{(t+1)} - x_1^{(t)} } + \left\|x_1^{(t+1)} - x_1^{(t)} \right\|^2.
        \end{aligned}
    \end{equation}
    For the inner product term, we have, by Equation \eqref{eq:lem:smooth-2},
    \begin{equation}
        \label{eq:lem:convex-2}
        \begin{aligned}
            &\ip{x_1^{(t)} - x^*}{x_1^{(t+1)} - x_1^{(t)}} \\
            = &  -\gamma \ip{x_1^{(t)} - x^*}{\bpi_\bR^\T \bPi_\bC \bY^{(t)}} - \gamma \ip{x_1^{(t)} - x^*}{\mone^\T \bY^{(t)}}.\\
        \end{aligned}
    \end{equation}
    Invoking Lemma \ref{lem:independ1} and Equation \eqref{eq:pp22}, we have $\expect \ip{x_1^{(t)} - x^*}{\bpi_\bR^\T \sum_{m=0}^{\min\crk{t,d_\bC}}\bA_m \bTh^{(t-m)}} = 0.$
    Then, denote $\bQ_{1}^{(t)} := \ip{x_1^{(t)} - x^*}{\bpi_\bR^\T \sum_{m=0}^{\min\crk{t,d_\bC}}\bA_m \nabla \bar{\bF}^{(t-m)} } ,$ and $\bQ_{2}^{(t)} := \ip{x_1^{(t)} - x^*}{\bpi_\bR^\T \sum_{m=0}^{\min\crk{t,d_\bC}}\bA_m \nabla \bF(\bar{\bX}^{(t-m)}) }.$
    It implies that
    \begin{equation}
        \label{eq:lem:convex-5}
        | \expect \bQ_{1}^{(t)}| \le \frac{ n \mu}{16} \norm{x_1^{(t)} - x^*}^2 + \frac{4}{n\mu} \expect \norm{\bpi_\bR^\T \sum_{m=0}^{\min\crk{t,d_\bC}}\bA_m \nabla \bar{\bF}^{(t-m)}}^2.
    \end{equation}
    Following similar procedures for getting Equation \eqref{eq:lem:smooth-4} and \eqref{eq:lem:xk-Ak-2-1}, we have
    \begin{equation}
        \label{eq:lem:convex-6}
        \begin{aligned}
             & -\gamma \expect \bQ_{2}^{(t)} \le   \frac{\gamma n \mu}{16} \norm{x_1^{(t)} - x^*}^2 - \gamma n \expect \ip{x_1^{(t)} - x^*}{\nabla f(x_1^{(t)})} \\
            & + \frac{4\gamma}{n\mu} \expect \norm{\bpi_\bR^\T \sum_{m=0}^{\min\crk{t,d_\bC}}\bA_m \nabla \bF(\bar{\bX}^{(t-m)}) - \mone^\T \nabla \bF(\bar{\bX}^{(t)})}^2. \\
        \end{aligned}
    \end{equation}
    Furthermore, from Equation \eqref{eq:lem:convex-2}, we have
    \begin{equation}
        \label{eq:lem:convex2-1}
        \begin{aligned}
            & -\gamma \expect \ip{x_1^{(t)} - x^*}{\mone^\T \bY^{(t)}} \\
            = & -\gamma \expect \ip{x_1^{(t)} - x^*}{\mone^\T \nabla \bar{\bF}^{(t)}} - \gamma n \expect \ip{x_1^{(t)} - x^*}{\nabla f(x_1^{(t)})}. \\
        \end{aligned}
    \end{equation}
    By Assumption \ref{a.convex},
    \begin{equation}
        \label{eq:lem:convex2-3}
        \begin{aligned}
        &\ip{x_1^{(t)} - x^*}{\nabla f(x_1^{(t)})} \ge \frac{\mu}{2}\norm{x_1^{(t)} - x^*}^2 + f(x_1^{(t)}) - f^*\\
        \ge & \frac{1}{2L}\norm{\nabla f(x_1^{(t)})}^2 + \frac{\mu}{2}\norm{x_1^{(t)} - x^*}^2.
        \end{aligned}
    \end{equation}
    Then, back to Equation \eqref{eq:lem:convex2-1}, it holds that
    \begin{equation}
        \label{eq:lem:convex2-4}
        \begin{aligned}
            & -\gamma \expect \ip{x_1^{(t)} - x^*}{\mone^\T \bY^{(t)}} \\
            & \le -\frac{\gamma n \mu}{4}\expect \norm{x_1^{(t)} - x^*}^2 + \gamma \kappa L\expect \norm{\hat{\bX}^{(t)}}_F^2 - \frac{\gamma n}{2L}\norm{\nabla f(x_1^{(t)})}^2
        \end{aligned} 
    \end{equation}
    From Equation \eqref{eq:lem:convex-1}, we have
    \begin{equation}
        \label{eq:lem:convex-7}
        \begin{aligned}
            & \expect \norm{x_1^{(t+1)} - x^*}^2 = \expect \norm{x_1^{(t)} - x^*}^2 - 2\gamma \expect\bQ_{1}^{(t)} - 2\gamma \expect\bQ_{2}^{(t)} \\
            & \qquad \qquad - 2\gamma \expect \ip{x_1^{(t)} - x^*}{\mone^\T \bY^{(t)}} + \expect \norm{x_1^{(t+1)}- x_1^{(t)}}_2^2 \\
            \le & \prt{1 - \frac{\gamma n \mu}{4}}  \expect \norm{x_1^{(t)} - x^*}^2 + \frac{8\gamma}{n\mu} \expect \norm{\bpi_\bR^\T \sum_{m=0}^{\min\crk{t,d_\bC}}\bA_m \nabla \bar{\bF}^{(t-m)}}^2 \\
            & + \frac{8\gamma}{n\mu} \expect \norm{\bpi_\bR^\T \sum_{m=0}^{\min\crk{t,d_\bC}}\bA_m \nabla \bF(\bar{\bX}^{(t-m)}) - \mone^\T \nabla \bF(\bar{\bX}^{(t)}) }^2 \\
            & + 2\gamma \kappa L\expect \norm{\hat{ \bX}^{(t)}}_F^2 - \frac{2 \gamma n}{L}\expect \norm{\nabla f(x_1^{(t)})}^2 + \frac{1}{n}\expect \norm{\Delta\bar{\bX}^{(t)}}_F^2.
        \end{aligned}
    \end{equation}
    Unwinding the recursion, 
    plugging Lemmas \ref{lem:Ak-barF} and \ref{lem:xk-Ak-2}, we have, for $\gamma \le 1/\prt{100 n\max\crk{d_\bR,d_\bC}\kappa L}$,
    \begin{equation}
        \label{eq:lem:convex-9}
        \begin{aligned}
            & \expect \norm{x_1^{(T)} - x^*}^2 \le \prt{1 - \frac{\gamma n \mu}{4}}^{T} \norm{x_1^{(0)} - x^*}^2 \\
            & + 8\gamma  \kappa L \sum_{t=0}^{T} \prt{1-\frac{\gamma n \mu}{4}}^{T-t}\brk{\prt{1-\frac{\gamma n \mu}{4}}^{-d_\bC}+ \frac{1}{4}} \expect\norm{\hat{\bX}^{(t)}}_F^2 \\
            & + \sum_{t=0}^{T} \prt{1-\frac{\gamma n \mu}{4}}^{T-t} \brk{50\expect\norm{x_1^{(t+1)} - x_1^{(t)} }^2- \frac{2\gamma n}{L}\expect \norm{\nabla f(x_1^{(t)})}^2}.\\
        \end{aligned}
    \end{equation}
    Plugging Lemma \ref{lem:X_diff} and \ref{lem:PiX} into Equation \eqref{eq:lem:convex-9} iteratively, we have, for $\gamma \le 1/\prt{1000 n\max\crk{d_\bR,d_\bC}\kappa L}$,
    \begin{equation}
        \label{eq:lem:convex-11}
        \begin{aligned}
            & \expect \norm{x_1^{(T)} - x^*}^2 \le \prt{1 - \frac{\gamma n \mu}{4}}^{T} \norm{x_1^{(0)} - x^*}^2 +  \frac{1200\gamma  \sigma^2}{ \mu}\\
            &  + 30000 \gamma^2 n \min\crk{d_\bR,d_\bC} r_\avg c_\avg \kappa^2 \sigma^2 + 800000 \gamma^3 n^2 d_\bR \kappa^2 L \sigma^2 \\
            &  + 800 \gamma^2 n \min\crk{d_\bR,d_\bC} r_\avg c_\avg \kappa^2 \sigma^2\prt{1-\frac{\gamma n \mu}{4}}^{-d_\bC} \\
            &  + 1800 \gamma^2 nc_\avg \kappa \prt{1-\frac{\gamma n \mu}{4}}^{T-d_\bC} \norm{ \nabla \bF^{(0)}}^2_F \\
            & + C_0\sum_{t=0}^{T} \prt{1-\frac{\gamma n \mu}{4}}^{T-t}\norm{\nabla f(x_1^{(t)})}^2 ,\\
        \end{aligned}
    \end{equation}
    where $C_0 := - \frac{2\gamma n}{L} + 900\gamma^2 n^2\kappa +  34 C \gamma^2 n^3 d_\bR \prt{1-\frac{\gamma n \mu}{4}}^{- d_\bC}  $ and $C :=8\gamma  \kappa L \prt{1-\frac{\gamma n \mu}{4}}^{-d_\bC} + 302\gamma \kappa L $.
    
    Then, for \smallmath{\gamma \le 1/\prt{1000 n\max\crk{d_\bR,d_\bC} r_\avg c_\avg\kappa L} },
    we have, by choosing $d = 20\max\crk{d_\bR,d_\bC} r_\avg c_\avg\kappa $ in Lemma \ref{lem:gamma}, that $\prt{1-\frac{\gamma n \mu}{4}}^{-d_\bC} \le \frac{1}{20\gamma n \max\crk{d_\bR,d_\bC} r_\avg c_\avg \kappa L}$, 
    and 
    \begin{equation}
        \label{eq:lem:convex-12}
        \begin{aligned}
            & 800000 \gamma^3 n^2 d_\bR \kappa^2 L \sigma^2 \le 1000 \gamma^2 n \min\crk{d_\bR,d_\bC} r_\avg c_\avg \kappa^2 \sigma^2,\\
            & 800 \gamma^2 n \min\crk{d_\bR,d_\bC} r_\avg c_\avg \kappa^2 \sigma^2\prt{1-\frac{\gamma n \mu}{4}}^{-d_\bC}  \le \frac{800 \gamma  \sigma^2}{\mu},\\
            & 20 \gamma^3 n^2 d_\bR  L \prt{1-\frac{\gamma n \mu}{4}}^{-d_\bC} \le  \gamma^2 n \min\crk{d_\bR,d_\bC} r_\avg c_\avg .\\
        \end{aligned}
    \end{equation}
    Besides, the coefficient of the term $\norm{ \nabla \bF^{(0)} }_F^2$ is bounded by $1800 \gamma^2 nc_\avg \kappa \prt{1-\frac{\gamma n \mu}{4}}^{-d_\bC} \le \frac{100\gamma}{L} \le \frac{1}{nL^2}$.
    
    Furthermore, for $\gamma \le \frac{1}{1000 n\max\crk{d_\bR,d_\bC} r_\avg c_\avg\kappa L}$, the coefficient of the term $\norm{\nabla f(x_1^{(t)})}^2$ is bounded by
    \begin{equation}
        \label{eq:lem:convex-14}
        \begin{aligned}
            & - \frac{2\gamma n}{L} + 900\gamma^2 n^2\kappa +  34 C \gamma^2 n^3 d_\bR \prt{1-\frac{\gamma n \mu}{4}}^{- d_\bC} \\
            \le & - \frac{2\gamma n}{L} + \frac{\gamma n}{L} + 320 \gamma^3 n^3 d_\bR \kappa L  \prt{1-\frac{\gamma n \mu}{4}}^{-2d_\bC} \\
            & + 12000 \gamma^3 n^3 d_\bR \kappa L  \prt{1-\frac{\gamma n \mu}{4}}^{- d_\bC} \\
            \le & - \frac{\gamma n}{L} + 320 \gamma^2 n^2  + 600 \gamma^2 n^2 \le - \frac{\gamma n}{L} + \frac{320\gamma n}{1000 L} +\frac{600\gamma n}{1000 L} \le 0.  
        \end{aligned}
    \end{equation}
    Then, back to Equation \eqref{eq:lem:convex-11}, we have, for $\gamma \le \frac{1}{1000 n\max\crk{d_\bR,d_\bC} r_\avg c_\avg\kappa L}$,
    \begin{equation}
        \label{eq:lem:convex-15}
        \begin{aligned}
            & \expect \norm{x_1^{(T)} - x^*}^2 \le \prt{1 - \frac{\gamma n \mu}{4}}^{T} \brk{\norm{x_1^{(0)} - x^*}^2 + \frac{1}{n L^2} \norm{\nabla \bF^{(0)} }_F^2 } \\
            & \quad +  \frac{2000\gamma  \sigma^2}{ \mu}  + 40000 \gamma^2 n \max\crk{d_\bR,d_\bC} r_\avg c_\avg \kappa^2 \sigma^2.
        \end{aligned}
    \end{equation}

\bibliography{references_all}
\bibliographystyle{siamplain}

\end{document}